\documentclass[11pt]{amsart}
\usepackage[T1]{fontenc}
\usepackage[utf8]{inputenc}
\usepackage{lmodern,amsmath,amssymb,mathtools,microtype,mathrsfs}
\usepackage[margin=29mm]{geometry}
\usepackage[hidelinks]{hyperref}
\newcommand{\R}{\mathbb R}\newcommand{\Z}{\mathbb Z}\newcommand{\N}{\mathbb N}\newcommand{\T}{\mathbb T}\newcommand{\C}{\mathbb C}\newcommand{\Q}{\mathbb Q}
\newcommand{\PP}{\mathbb P}\newcommand{\e}{\mathrm e}\newcommand{\dd}{\,\mathrm d}
\newcommand{\PM}{\operatorname{PM}}\newcommand{\PF}{\operatorname{PF}}\newcommand{\supp}{\operatorname{supp}}\newcommand{\rank}{\operatorname{rank}}\newcommand{\tr}{\operatorname{tr}}\newcommand{\dist}{\operatorname{dist}}\newcommand{\calF}{\mathcal F}\newcommand{\HH}{\mathcal H}\newcommand{\tn}[1]{\|#1\|_{\T}}\newcommand{\pr}[1]{\langle#1\rangle_0}\newcommand{\shift}{\mathsf T}
\newtheorem{theorem}{Theorem}[section]\newtheorem{lemma}[theorem]{Lemma}\newtheorem{proposition}[theorem]{Proposition}\newtheorem{corollary}[theorem]{Corollary}
\theoremstyle{definition}\newtheorem{definition}[theorem]{Definition}
\theoremstyle{remark}\newtheorem{remark}[theorem]{Remark}
\numberwithin{equation}{section}
\title[Polars of sets of uniqueness]{Separating the polar of sets of uniqueness\protect\\from the polars of all $H^{(d)}$-sets}
\author{V\'aclav Vlas\'ak}
\thanks{Corresponding author: V\'aclav Vlas\'ak. Faculty of Mathematics and Physics, Charles University, Sokolovsk\'a 83, 186 75 Praha 8, Czech Republic. E-mail: \href{mailto:vaclav.vlasak@matfyz.cuni.cz}{\texttt{vaclav.vlasak@matfyz.cuni.cz}}.}
\date{November 19, 2026}
\hypersetup{pdftitle={Separating the polar of sets of uniqueness from the polars of all H(d)-sets},pdfauthor={Václav Vlasák}}
\begin{document}
\begin{abstract}
We construct a perfect set of uniqueness on the circle carrying a continuous probability measure which annihilates every $H^{(d)}$-set, simultaneously for all finite dimensions $d$. This gives a strict inclusion between the polar of the compact sets of uniqueness and the polar of the union of the classes $H^{(d)}$, answering a question raised in earlier work on these polars. The construction uses random binary increments of size comparable to the reciprocal of a product of square roots. A Gaussian sampling estimate gives a uniform Fourier-tail bound for all pseudomeasures supported by the resulting set, and also bounded atomic approximation and spectral synthesis. A separate energy estimate, stable under deletion of a small number of increments, yields full support for every limiting quasi-independent image of every positive measure absolutely continuous with respect to the constructed measure.
\end{abstract}
\maketitle

\section{Introduction}
The classes of sets of uniqueness and of extended uniqueness describe different forms of thinness for trigonometric series. Their polars provide a measure-theoretic way of comparing them. For a family $\mathcal A$ of compact subsets of $\T=\R/\Z$, write
\[
\mathcal A^\perp=\{\nu\in M_+(\T):\nu(K)=0\text{ for every }K\in\mathcal A\},
\]
where $\N=\{1,2,\ldots\}$ and $M_+(\T)$ denotes the finite positive Borel measures. We denote by $U$ the compact sets of uniqueness and by $U_0$ the compact sets of extended uniqueness. A compact set belongs to $U$ if a trigonometric series whose symmetric partial sums tend to zero at every point of its complement has all coefficients zero. The class $U_0$ consists of the compact sets annihilated by every Rajchman measure, that is, every finite measure whose Fourier coefficients tend to zero at infinity.

Lyons identified $U_0^\perp$ with the positive Rajchman measures \cite{Lyons84} and proved that the polar of the union of the classes $H^{(d)}$ strictly contains this class \cite{Lyons87}. Kaufman proved that $U^\perp\ne U_0^\perp$ \cite{Kaufman}. These results do not decide whether replacing $U_0$ by $U$ still separates it from the union of all $H^{(d)}$.

To specify the latter classes, put $\e(x)=\exp(2\pi ix)$ for $x\in\R$ and $\T^d=\R^d/\Z^d$. For $a\in\Z^d$, let $m_a:\T\to\T^d$ be the continuous homomorphism $m_a(x)=ax\pmod{\Z^d}$. We use the terminology of \cite{V}.
\begin{samepage}
\begin{definition}\label{def:H}
A sequence $(a^j)_{j\ge1}$ in $\Z^d$ is quasi-independent if
\[
\forall\alpha\in\Z^d\setminus\{0\},\qquad |\alpha\cdot a^j|\longrightarrow\infty.
\]
For $d\ge1$, $H^{(d)}$ is the family of compact sets $K\subset\T$ for which there are a quasi-independent sequence $(a^j)$ in $\N^d$ and a nonempty open rectangle $J\subset\T^d$ such that
\begin{equation}\label{eq:gap}
m_{a^j}(K)\cap J=\varnothing\qquad(j\ge1).
\end{equation}
Set $\HH=\bigcup_{d\ge1}H^{(d)}$.
\end{definition}
\end{samepage}
The full-support theorem below allows arbitrary integer vectors, so in particular applies to this convention for $H^{(d)}$. The classical inclusions are $H^{(d)}\subset H^{(d+1)}\subset U\subset U_0$; see \cite{KL,V}. It was shown in \cite{V} that $(H^{(d+1)})^\perp\subsetneq(H^{(d)})^\perp$ for every $d$. In particular, each individual $(H^{(d)})^\perp$ differs from $U^\perp$. The simultaneous question, explicitly raised in \cite[p.~202]{V}, is settled by the following theorem.

\begin{theorem}\label{thm:main}
There exist a perfect compact set $E\subset\T$ and a continuous probability measure $\mu$ with $\supp\mu=E$ such that
\[
E\in U,\qquad \mu(K)=0\quad(K\in H^{(d)},\ d\ge1).
\]
Consequently $U^\perp\subsetneq\HH^\perp$. Moreover, $E$ can be chosen to be a set of spectral synthesis, with uniformly bounded finite-rank atomic approximation of its supported pseudomeasures.
\end{theorem}
Here continuity of a measure means absence of atoms. The precise approximation statement, including its norm bound, is given in Proposition~\ref{prop:synthesis}.

The construction uses two product spaces, with distinct roles. The coefficient space is
\[
\Omega=[1,2]^{\N},\qquad
\PP_{\rm coef}=\bigotimes_{k\ge1}\bigl(\mathcal L^1|_{[1,2]}\bigr),
\]
with its product Borel sigma-algebra; $\mathcal L^1|_{[1,2]}$ has total mass one. Thus the coordinate functions $\omega\mapsto\omega_k$ are independent and uniform on $[1,2]$. Set
\begin{equation}\label{eq:model}
M_k=\lceil16\sqrt{k}\rceil,\qquad Q_0=1,\qquad
Q_n=\prod_{k=1}^nM_k,\qquad
 t_k(\omega)=\frac{\omega_k}{Q_k}\quad(k\ge1).
\end{equation}
Once $\omega\in\Omega$ is fixed, we abbreviate $t_k(\omega)$ to $t_k$. The digit space is $\mathcal Z=\{0,1\}^{\N}$, with fair product probability
$\PP_{\rm dig}=\bigotimes_{k\ge1}(\frac12\delta_0+\frac12\delta_1)$ and coordinate functions $Z_k(z)=z_k$.
We equip topological spaces with their Borel sigma-algebras. For a Borel map $p:X\to Y$ and a Borel measure $\nu$ on $X$, the notation $p_*\nu$ means the image Borel measure:
$(p_*\nu)(B)=\nu(p^{-1}(B))$ for every Borel set $B\subset Y$.
For each fixed coefficient sequence define
\begin{equation}\label{eq:coding}
\begin{gathered}
\pi_\omega:\mathcal Z\longrightarrow\R,\qquad
\pi_\omega(z)=\sum_{k\ge1}t_k(\omega)z_k,\\
E_\omega=\pi_\omega(\mathcal Z),\qquad
\mu_\omega=(\pi_\omega)_*\PP_{\rm dig}.
\end{gathered}
\end{equation}
We omit the subscript $\omega$ when the coefficients are fixed. Lemma~\ref{lem:geometry} shows that $E_\omega\subset[0,2/15]$ and that $\pi_\omega$ is a homeomorphism onto its image. We identify this interval with its image in $\T$; the same notation $E_\omega,\mu_\omega$ is used for the resulting set and measure on the circle. The conclusions of Theorem~\ref{thm:main} hold for $\PP_{\rm coef}$-almost every $\omega\in\Omega$. The randomness of the coefficients is used to select $\omega$; after that selection, the digits alone describe the measure $\mu_\omega$.

The proof has two parts after the coefficients are fixed. For uniqueness and synthesis, we localize a supported pseudomeasure to the binary pieces of $E$ and replace each piece by a Dirac measure at its base point. The coefficients of these atomic approximations need not be positive; the required uniform bound is in the pseudomeasure norm. The auxiliary polynomials in Sections~\ref{sec:sampling}--\ref{sec:PM} encode their Fourier transforms, while differentiation in one real parameter controls the approximation error. Bounded atomic approximation yields spectral synthesis. Uniqueness follows separately from the Fourier-tail estimate below, not from atomic approximability alone. A small dual Gaussian sum gives a dimension-independent sampling estimate for binary trigonometric polynomials. Localization then gives, on the line, the estimate below. Here $\PM(\R)$ denotes the space of pseudomeasures, with Fourier-transform norm as defined in Section~\ref{sec:tools}:
\begin{equation}\label{eq:tailgoal}
\|S\|_{\PM(\R)}\le80\limsup_{q\to+\infty}|\widehat S(q)|\qquad(\supp S\subset E).
\end{equation}
This implies uniqueness and bounded synthesis. A second estimate prevents the arithmetic energy of the coefficients from being concentrated on a small number of coordinates. In an induction on dimension it produces a nonzero Gaussian factor whenever the finite-energy subgroup has full rank. Its support supplies the additional invariance needed to pass to a lower-dimensional quotient. This proves a full-support assertion for every positive $\nu\ll\mu$, from which the polar conclusion follows.

The synthesis assertion should be compared with classical results for symmetric sets. Schneider's Theorem~1.16 \cite{Schneider} assumes square summability of the successive increment ratios. Here
\[
\frac1{2M_{k+1}}\le\frac{t_{k+1}}{t_k}\le\frac2{M_{k+1}},\qquad
\sum_k(t_{k+1}/t_k)^2=\infty.
\]
That criterion therefore does not provide the estimates used here. The principal conclusion is the simultaneous polar separation; synthesis is an additional property of the same construction.

\section{Conventions and standard tools}\label{sec:tools}
\subsection{Fourier duality and localization}
For a vector in $\R^d$ or $\C^d$, $|\cdot|$ denotes the Euclidean norm. For a scalar function $f$ on a specified set $X$, we use $\|f\|_{\infty,X}=\sup_{x\in X}|f(x)|$; a shorter $\|f\|_\infty$ will be used only after specifying $X$. The notation $\|f\|_{L^p(\nu)}$ denotes the norm with respect to the stated measure $\nu$. Lebesgue measure on $\R$ is understood when the argument of a Fourier-transform norm is a function on $\R$.
We use the Fourier-duality framework of \cite{KL}, writing $\PM$ for pseudomeasures and $\PF$ for pseudofunctions. All dual pairings below are complex-linear in the test function. Write
\[
A(\T)=\left\{f=\sum_{k\in\Z}a_k\e(kx): (a_k)\in\ell^1(\Z)\right\},\qquad \|f\|_A=\sum_k|a_k|.
\]
Thus $\PM(\T)=A(\T)^*\simeq\ell^\infty(\Z)$, with complex-linear pairing
\[
\widehat S(k)=\langle S,\e(k\,\cdot)\rangle,\qquad
\langle S,f\rangle=\sum_ka_k\widehat S(k).
\]
In particular,
\[
\|S\|_{\PM(\T)}=\sup_{k\in\Z}|\widehat S(k)|.
\]
The series defining $f\in A(\T)$ and its pairing with $S$ are absolutely convergent. The subspace $\PF(\T)$ consists of those $S$ for which $\widehat S(k)\to0$ as $|k|\to\infty$, and thus corresponds to $c_0(\Z)$. These identifications concern Fourier coefficients, not functions on the original circle. The classical pseudofunction characterization is
\begin{equation}\label{eq:U}
E\in U\quad\Longleftrightarrow\quad \PM(E)\cap\PF(\T)=\{0\},
\end{equation}
where $\PM(E)=\{S\in\PM(\T):\supp S\subset E\}$; see \cite[Chapter II]{KL} and \cite[Section 2.1]{MZ}.

For a finite measure $\nu$ on $\T^d$ we use the same sign convention, $\widehat\nu(\alpha)=\int_{\T^d}\e(\alpha\cdot x)\dd\nu(x)$ for $\alpha\in\Z^d$. In particular, the plus sign in the coefficient is intentional.

On $\R$ we use $A(\R)=\{\int a(q)\e(qx)\dd q:a\in L^1(\R)\}$ with norm $\|a\|_1$, and $\PM(\R)=A(\R)^*\simeq L^\infty(\R)$. The Fourier transform has the plus-sign convention
$\widehat f(q)=\calF f(q)=\int f(x)\e(qx)\dd x$; it extends by duality to tempered distributions. We write $\widehat S=\calF S$, and use $\calF$ for longer expressions. Pseudofunctions on $\R$ correspond to transforms with representatives in $C_0(\R)$, the continuous functions vanishing at infinity. The norm is
\[
\|S\|_{\PM(\R)}=\|\widehat S\|_{L^\infty(\R)}.
\]
Initially this is an essential supremum. When $S$ has compact support, its transform has a continuous representative, as recalled below, and the norm equals $\sup_{q\in\R}|\widehat S(q)|$. We always use that representative for pointwise evaluations.

Smooth means $C^\infty$. Distributional test functions on $\R$ belong to $\mathcal D=C_c^\infty$; tempered distributions act on the Schwartz space $\mathcal S$. On the torus the test functions are $C^\infty(\T^d)$. Supports of pseudomeasures are their distributional supports. The transform of a compactly supported distribution has a smooth (indeed entire) representative \cite[Theorem 7.23(a)]{RudinFA}. In particular, for a cutoff $\chi=1$ near $\supp S$,
\[
\widehat S(q)=\langle S,\chi(x)\e(qx)\rangle.
\]
The standard product--convolution identity and its norm estimate are
\begin{equation}\label{eq:multiplier}
\widehat{\phi S}=\widehat\phi*\widehat S,\qquad
\|\phi S\|_{\PM(\R)}\le\|\widehat\phi\|_{L^1(\R)}\|S\|_{\PM(\R)}\quad(\phi\in C_c^\infty(\R)).
\end{equation}
See \cite[Theorems 7.15 and 7.19]{RudinFA}, with normalization adjusted as above. On the torus the corresponding rule is the $\ell^1*\ell^\infty$ estimate; $\ell^1*c_0\subset c_0$ also gives localization of pseudofunctions.

\subsection{Measures, convergence, and Gaussian measures}\label{sec:prob}
For topological spaces $X,Y$, a finite positive Borel measure $\nu$ on $X$, and a Borel map $p:X\to Y$, the image measure defined in the introduction satisfies
\[
\int_Y\varphi\dd(p_*\nu)=\int_X\varphi(p(x))\dd\nu(x)
\]
for every bounded Borel function $\varphi:Y\to\C$. If $g:X\to[0,\infty)$ is a Borel function integrable with respect to $\nu$, $g\nu$ denotes the measure $B\mapsto\int_Bg\dd\nu$; thus $p_*(g\nu)$ has mass $\int_Xg\dd\nu$.
On a probability space $(X,\mathscr A,\PP)$, a map $X_0:X\to\R^d$ measurable from $\mathscr A$ to the Borel sigma-algebra of $\R^d$ is a random vector. Probability theory calls $(X_0)_*\PP$ its \emph{law}. We will use image-measure notation, also when a weight is present.

On $\T^d$, convolution is the image of a product measure under the continuous addition map $a:\T^d\times\T^d\to\T^d$, $a(v,y)=v+y$:
$\kappa*\gamma=a_*(\kappa\otimes\gamma)$. In particular, for $x,y\in\T^d$, we have $\delta_x*\delta_y=\delta_{x+y}$, where addition is understood in $\T^d$. Consequently, if $V,Y$ are torus-valued random vectors, $g\ge0$ is integrable, and the pair $(V,g)$ is independent of $Y$, then
\begin{equation}\label{eq:weightedlaw}
(V+Y)_*(g\PP)=V_*(g\PP)*Y_*\PP.
\end{equation}
Indeed, for $\varphi\in C(\T^d)$, independence gives
\[
\int g\,\varphi(V+Y)\dd\PP
=\int\!\int\varphi(v+y)\dd[V_*(g\PP)](v)\dd[Y_*\PP](y),
\]
which is the integral against the right-hand side of \eqref{eq:weightedlaw}.

We write $\nu_j\xrightarrow{w^*}\nu$ in $M(\T^d)=C(\T^d)^*$ to mean
\[
\int_{\T^d}\varphi\dd\nu_j\longrightarrow\int_{\T^d}\varphi\dd\nu
\qquad\text{for every }\varphi\in C(\T^d).
\]
Here $M(\T^d)$ is the space of finite complex Borel measures with the total variation norm. For a positive measure that norm is its mass $\nu(\T^d)$. For uniformly bounded total variation, testing all characters is sufficient: first use linearity for trigonometric polynomials, and then their uniform density in $C(\T^d)$. Positive measures of uniformly bounded mass have weak-star convergent subsequences by Banach--Alaoglu and separability of $C(\T^d)$. Their limits are positive, and testing the constant function one gives convergence of the masses.

Let $\rho_d:\R^d\to\T^d$, $\rho_d(x)=x+\Z^d$, be the canonical projection. For a real symmetric positive semidefinite matrix $C$, let $C^{1/2}$ denote its unique symmetric positive semidefinite square root.
The possibly singular Gaussian measure $N(0,C)$ is the image of the standard Gaussian measure $N(0,I_d)$ under $x\mapsto C^{1/2}x$. Its covariance matrix is $C$. For $\gamma_C=(\rho_d)_*N(0,C)$, the Gaussian Fourier transform gives
\begin{equation}\label{eq:gaussianFourier}
\widehat\gamma_C(\alpha)=\exp(-2\pi^2\alpha^{\mathsf T}C\alpha)
\qquad(\alpha\in\Z^d),
\end{equation}
including when $C$ is singular.

We will use Stirling's formula in the single form
\begin{equation}\label{eq:Stirling}
\log\Gamma(x+1)=x\log x-x+\tfrac12\log(2\pi x)+O(x^{-1})\qquad(x\to\infty),
\end{equation}
see \cite[5.11.1]{DLMF}.

\section{Geometry and two almost-sure estimates}\label{sec:random}
Fix $\omega\in\Omega$, so that $\omega_k\in[1,2]$ and $t_k=\omega_k/Q_k$ as in \eqref{eq:model}. For the map $\pi=\pi_\omega$ and set $E=E_\omega$ in \eqref{eq:coding}, put
\[
L_n=\sum_{k>n}t_k,\qquad x_z=\sum_{k\le n}t_kz_k,\qquad
E_z=\pi([z])\quad(z\in\{0,1\}^n),
\]
where $[z]=\{w\in\mathcal Z:w_k=z_k\ (1\le k\le n)\}$ is a cylinder in the digit space. Thus $z\in\{0,1\}^n$ is a finite binary vector, $[z]\subset\mathcal Z$ is a cylinder, and $E_z\subset E$ is its image. For any infinite sequence $w$, write $w|_n=(w_1,\ldots,w_n)$ for its prefix of length $n$.
\begin{lemma}\label{lem:geometry}
For every coefficient sequence in $[1,2]^{\N}$,
\begin{gather}
E\subset[0,2/15],\qquad L_n\le\frac8{3Q_nM_{n+1}}\le\frac1{6Q_n\sqrt{n+1}},\qquad t_n>L_n,\label{eq:tails}\\
\log Q_n=n\log16+\tfrac12n\log n-\tfrac12n+o(n).\label{eq:growth}
\end{gather}
The map $\pi$ is a homeomorphism onto $E$, and $\mu(E_z)=2^{-n}$. In particular $E$ is perfect and $\mu$ is continuous with support $E$. For each $n$, the sets $E_z$ lie in the intervals $[x_z,x_z+L_n]$, which are pairwise separated by positive distances.
\end{lemma}
\begin{proof}
Since $M_k\ge16$ and $t_k\le2/Q_k$, summation of a geometric series gives $\sum t_k\le2/15$ and
$L_n\le32/(15Q_{n+1})\le8/(3Q_nM_{n+1})$. The remaining tail assertions follow from $M_{n+1}\ge16\sqrt{n+1}$ and $t_n\ge1/Q_n$.
The series defining $\pi$ converges uniformly. If $z,w$ first differ at $m=\min\{k:z_k\ne w_k\}$, then
$|\pi(z)-\pi(w)|\ge t_m-L_m>0$. Thus $\pi$ is a continuous injection from a compact space and is therefore a homeomorphism onto its image. The product measure gives every cylinder positive mass and every singleton zero mass. The asserted properties follow. For two finite binary vectors of length $n$ whose first differing coordinate is $m$, the distance between the corresponding intervals is at least $t_m-\sum_{m<k\le n}t_k-L_n=t_m-L_m>0$.
To check the growth estimate, first write $M_k=16\sqrt{k}(1+r_k)$, where $0\le r_k<(16\sqrt{k})^{-1}$. The inequality $0\le\log(1+r)\le r$ gives
\begin{equation}\label{eq:logMk}
\log M_k=\log16+\tfrac12\log k+O(k^{-1/2}).
\end{equation}
Consequently,
\begin{align*}
\log Q_n
&=\log\left(\prod_{k=1}^nM_k\right)=\sum_{k=1}^n\log M_k\\
&\stackrel{\eqref{eq:logMk}}{=}n\log16+\tfrac12\sum_{k=1}^n\log k
   +O\left(\sum_{k=1}^nk^{-1/2}\right)\\
&=n\log16+\tfrac12\log(n!)+O(\sqrt n)\\
&\stackrel{\eqref{eq:Stirling}}{=}n\log16+\tfrac12n\log n-\tfrac12n
  +\tfrac14\log(2\pi n)+O(\sqrt n).
\end{align*}
Since $\log n=o(n)$ and $\sqrt n=o(n)$, this proves \eqref{eq:growth}.
\end{proof}

\subsection{The dual Gaussian sum}
We keep the fixed smoothing constant $\varepsilon=1/4$ explicit in the Gaussian formulas, so that the normalization of the Gaussian integrals remains visible. It is not an additional free parameter. For $t\in\R^n$ and $T>0$ define
\begin{equation}\label{eq:Delta}
D_n(t,T)=\sum_{b\in(\frac12\Z)^n\setminus\{0\}}
\exp\left(-\frac{\pi^2\varepsilon^2}{n}|b|^2-\pi^2T^2(b\cdot t)^2\right).
\end{equation}
The scalar product in this formula is real, without reduction modulo one. Put
\begin{equation}\label{eq:AB}
A_n=\sum_{r\in\frac12\Z}e^{-\pi^2\varepsilon^2r^2/n},\qquad
B_n=\sum_{r\in\frac12\Z\setminus\{0\}}\frac{e^{-\pi^2\varepsilon^2r^2/n}}{|r|}.
\end{equation}
All summands are nonnegative, and Gaussian decay gives $A_n<\infty$. Dropping the second term in the exponent and factoring the product series gives
\[
D_n(t,T)\le
\sum_{b\in(\frac12\Z)^n\setminus\{0\}}e^{-\pi^2\varepsilon^2|b|^2/n}
=\prod_{k=1}^n\left(\sum_{b_k\in\frac12\Z}e^{-\pi^2\varepsilon^2b_k^2/n}\right)-1
=A_n^n-1<\infty.
\]
This also proves absolute convergence of the lattice sum in \eqref{eq:Delta} for every $t\in\R^n$ and $T>0$. For a fixed deterministic $\tau>0$ write
\begin{equation}\label{eq:T}
T_n(\tau)=\tau Q_n\sqrt n,\qquad
\Delta_n(v,\tau)=D_n\left((v_k/Q_k)_{k=1}^n,T_n(\tau)\right),
\quad v\in[1,2]^n.
\end{equation}
\begin{lemma}\label{lem:Delta}
There is a constant $C$, independent of $n\ge1$ and $\tau>0$, such that
\[
\int_\Omega\Delta_n(\omega|_n,\tau)\dd\PP_{\rm coef}(\omega)\le\frac{C(1+\log n)}{\tau\sqrt n}.
\]
For each fixed $\tau>0$, the event
\begin{equation}\label{eq:D}
\mathcal D_\tau=\left\{\omega\in\Omega:\lim_{j\to\infty}\Delta_{2^j}(\omega|_{2^j},\tau)=0\right\}
\end{equation}
has $\PP_{\rm coef}$-measure one.
\end{lemma}
\begin{proof}
Fix $n,\tau$ and abbreviate $T=T_n(\tau)$. For $b\in(\frac12\Z)^n\setminus\{0\}$ let $m=\min\{k:b_k\ne0\}$. Hold the coordinates $\omega_k$, $k\ne m$, fixed and put $c_m=\sum_{k\ne m}b_k\omega_k/Q_k$. Integration in the remaining coordinate gives
\begin{equation}\label{eq:conditional}
\begin{aligned}
\int_1^2\exp\left[-\pi^2T^2\left(\frac{b_m v}{Q_m}+c_m\right)^2\right]\dd v
&=\frac{Q_m}{|b_m|}\int_{J_m}e^{-\pi^2T^2y^2}\dd y\\
&\le\frac{Q_m}{|b_m|}\int_\R e^{-\pi^2T^2y^2}\dd y
=\frac{Q_m}{\sqrt\pi T|b_m|},
\end{aligned}
\end{equation}
where $J_m$ is the interval which is the image of $[1,2]$ under $v\mapsto b_mv/Q_m+c_m$. The bound is independent of the coordinates that were held fixed. In particular, the translation $c_m$ disappears only when $J_m$ is enlarged to $\R$.
Put $a=\pi^2\varepsilon^2/4$. Writing the positive half-integers as $k/2$ and using that $x\mapsto e^{-ax^2/n}$ is decreasing on $[0,\infty)$, we obtain
\begin{equation}\label{eq:Abound}
\begin{aligned}
A_n
&=1+2\sum_{k=1}^{\infty}e^{-ak^2/n}\\
&\le1+2\int_0^\infty e^{-ax^2/n}\dd x\\
&=1+\sqrt{\frac{\pi n}{a}}
 =1+\frac{2\sqrt n}{\sqrt\pi\varepsilon}\\
&\le 8\sqrt n.
\end{aligned}
\end{equation}
The last inequality uses $\varepsilon=1/4$ and $n\ge1$.
For the second sum, put $K=\lceil\sqrt n\rceil$. The function $x\mapsto e^{-ax^2/n}/x$ is decreasing on $(0,\infty)$, so splitting the sum at $K$ gives
\begin{equation}\label{eq:Bbound}
\begin{aligned}
B_n
&=4\sum_{k=1}^{\infty}\frac{e^{-ak^2/n}}k\\
&=4\sum_{k=1}^{K}\frac{e^{-ak^2/n}}k
 +4\sum_{k=K+1}^{\infty}\frac{e^{-ak^2/n}}k\\
&\le4\sum_{k=1}^{K}\frac1k
 +4\int_K^\infty\frac{e^{-ax^2/n}}x\dd x\\
&\le4\left(1+\int_1^K\frac{\dd x}{x}
 +\int_{K/\sqrt n}^\infty\frac{e^{-av^2}}v\dd v\right)\\
&\le4\left(1+\log K+\int_1^\infty\frac{e^{-av^2}}v\dd v\right)\\
&\le4\left(1+\log2+\tfrac12\log n
 +\int_1^\infty\frac{e^{-av^2}}v\dd v\right)\\
&\le C'(1+\log n).
\end{aligned}
\end{equation}
Here the change of variables is $v=x/\sqrt n$, and $\sqrt n\le K\le2\sqrt n$. The remaining integral is finite and independent of $n$; hence $C'$ is a constant independent of $n$.
Monotone convergence and grouping by $m$ in \eqref{eq:conditional} yield
\begin{align}
\int_\Omega\Delta_n(\omega|_n,\tau)\dd\PP_{\rm coef}(\omega)
&\le\frac1{\sqrt\pi T}\sum_{m=1}^nQ_m
\sum_{b_m\in\frac12\Z\setminus\{0\}}\frac{e^{-\pi^2\varepsilon^2b_m^2/n}}{|b_m|}
\prod_{k=m+1}^n\sum_{b_k\in\frac12\Z}e^{-\pi^2\varepsilon^2b_k^2/n}\notag\\
&=\frac{B_n}{\sqrt\pi T}\sum_{m=1}^n Q_m A_n^{n-m}.\label{eq:DeltaSum}
\end{align}
Put $r=n-m$. For $r\ge1$, the geometric mean of $m+1,\ldots,n$ is at least $m+1$. Adjoining the smaller numbers $1,\ldots,m$ therefore cannot increase it. Hence
\[
\left(\frac{n!}{m!}\right)^{1/r}\ge(n!)^{1/n}\ge n/e,
\]
where the last step uses the standard bound $(n!)^{1/n}\ge n/e$.
Thus $n!/(n-r)!\ge(n/e)^r$. For $r=0$, this also holds, with the empty product equal to one. Consequently
\begin{equation}\label{eq:theta}
\frac{Q_m A_n^{n-m}}{Q_n}
\le\left(\frac{8\sqrt n}{16}\right)^r\sqrt{\frac{(n-r)!}{n!}}
\le\left(\frac{\sqrt e}{2}\right)^r.
\end{equation}
Since $\sqrt e/2<1$, combining \eqref{eq:DeltaSum}, \eqref{eq:Bbound}, and \eqref{eq:theta} now gives the complete estimate
\begin{align*}
\int_\Omega\Delta_n(\omega|_n,\tau)\dd\PP_{\rm coef}
&\le\frac{B_nQ_n}{\sqrt\pi T_n(\tau)}\sum_{r=0}^{n-1}\left(\frac{\sqrt e}{2}\right)^r\\
&\le\frac{C'(1+\log n)}{\sqrt\pi\tau\sqrt n}\frac1{1-\sqrt e/2}.
\end{align*}
For the dyadic indices, monotone convergence therefore yields
\begin{align*}
\int_\Omega\sum_{j\ge0}\Delta_{2^j}(\omega|_{2^j},\tau)\dd\PP_{\rm coef}
&=\sum_{j\ge0}\int_\Omega\Delta_{2^j}(\omega|_{2^j},\tau)\dd\PP_{\rm coef}\\
&\le\frac C\tau\sum_{j\ge0}(1+j\log2)2^{-j/2}<\infty.
\end{align*}
A nonnegative function with finite integral is finite almost everywhere. Hence the series inside the integral is finite for almost every $\omega$, and its terms tend to zero as $j\to\infty$. This proves $\PP_{\rm coef}(\mathcal D_\tau)=1$ for the fixed $\tau$.
\end{proof}

\subsection{Energy after deletion of coordinates}
For $x\in\R$ put $\tn{x}=\dist(x,\Z)$. The following estimate says that a fixed amount of arithmetic energy remains even after deleting any small set of coordinates.
\begin{lemma}\label{lem:energy}
There is a Borel event $\mathcal R\subset\Omega$ with $\PP_{\rm coef}(\mathcal R)=1$ such that for every $\omega\in\mathcal R$ there exists $n_*\in\N$ with the following property: for every $n\ge n_*$, every $k\in\Z$ with $Q_n\le|k|<Q_{n+1}$, and every $I\subset\{1,\ldots,n\}$ with $|I|\le\lfloor n^{3/4}\rfloor$,
\begin{equation}\label{eq:robust}
\sum_{\substack{\ell\le n\\\ell\notin I}}\tn{k\omega_\ell/Q_\ell}^2\ge c_{\mathrm{en}},
\qquad c_{\mathrm{en}}=\frac{2^{-14}}{\pi}.
\end{equation}

\end{lemma}
\begin{proof}
For fixed $n,k$ in this band, $k\omega_\ell/Q_\ell$ is uniform on an interval of length $|k|/Q_\ell\ge1$ for $\ell\le n$. Its reduction modulo one has density at most $2$, and these reductions remain independent. Put $s=\lfloor n^{3/4}\rfloor$, $m=n-s$. After deleting a fixed $s$ coordinates, the probability that the remaining squared distances sum to less than $c_{\mathrm{en}}$ is at most $2^mV_mc_{\mathrm{en}}^{m/2}$, where $V_m=\pi^{m/2}/\Gamma(m/2+1)$ is the volume of the unit Euclidean ball in $\R^m$. A failure after fewer deletions remains a failure after enlarging the deleted set. The union bound over all deleted sets and all $k$ therefore gives
\begin{equation}\label{eq:pn}
p_n=2Q_{n+1}\binom ns2^mV_mc_{\mathrm{en}}^{m/2}
\end{equation}
as an upper bound for a failure anywhere in the band.
Since $s\le n^{3/4}$, we have $s\log n=o(n)$, $m/n\to1$ and $m\log m=n\log n+o(n)$. Formula \eqref{eq:Stirling} gives
\[
\log V_m=\frac m2\log\pi-\log\Gamma(m/2+1)
=\frac m2\log(2\pi e/m)+O(\log m),
\]
while $\log\binom ns\le s\log n=o(n)$. Substituting these estimates into \eqref{eq:pn}, and using \eqref{eq:growth} also for $Q_{n+1}$, gives
\begin{align*}
\log p_n
&=\log2+\log Q_{n+1}+\log\binom ns
   +m\log2+\log V_m+\tfrac m2\log c_{\mathrm{en}}\\
&=n\log16+\tfrac12n\log n-\tfrac12n
  +n\log2+\tfrac n2\log(2\pi e)-\tfrac12n\log n
  +\tfrac n2\log c_{\mathrm{en}}+o(n).
\end{align*}
The $n\log n$ terms cancel. Hence
\[
\frac{\log p_n}{n}\longrightarrow\log16+\tfrac12\log(8\pi c_{\mathrm{en}})=-\tfrac12\log8<0.
\]
The bounds are summable, and the first Borel--Cantelli lemma proves the assertion. No independence between the failure events is used.
\end{proof}

\section{A Gaussian sampling estimate}\label{sec:sampling}
For the base points $x_z=z\cdot t$ from Section~\ref{sec:random}, a binary trigonometric polynomial and an atomic measure are related by
\[
P(\theta)=\sum_{z\in\{0,1\}^n}a_z\e(z\cdot\theta),\qquad
\lambda=\sum_{z\in\{0,1\}^n}a_z\delta_{x_z},\qquad
P(qt)=\sum_z a_z\e(qx_z)=\widehat\lambda(q).
\]
The next lemma controls the supremum norm of a binary trigonometric polynomial by its values on a short segment of the line in direction $t$. In Section~\ref{sec:PM}, we apply this estimate to polynomials arising from localized pseudomeasures to obtain the norm bounds needed for approximation by atomic measures. Its localization near an arbitrary prescribed $u$ is essential: in Lemma~\ref{lem:M}, the error of changing the second parameter of $P_n$ is proportional to $|q-u|$. A large value somewhere on the line, with no control of that distance, would not suffice. We state the estimate for arbitrary $t\in\R^n$; for the construction it will be applied to $t=(t_1,\ldots,t_n)$.
\begin{lemma}\label{lem:sampling}
Let $t\in\R^n$, $T>0$, and suppose
\begin{equation}\label{eq:admissible}
D_n(t,T)\le\tfrac12e^{-\pi^2\varepsilon^2}.
\end{equation}
For every $P:\R^n\to\C$ of the form $P(\theta)=\sum_{z\in\{0,1\}^n}a_z\e(z\cdot\theta)$ and every $u\in\R$, there is $q\in\R$ such that
\begin{equation}\label{eq:sample}
|q-u|\le20T,\qquad |P(qt)|\ge\|P\|_\infty/20.
\end{equation}
\end{lemma}
Here $\|P\|_\infty=\sup_{\theta\in\R^n}|P(\theta)|$, equivalently the supremum on $\T^n$. To interpret the conclusion geometrically, consider the segment $\{qt:|q-u|\le20T\}$ of the line in direction $t$. It has Euclidean length $40T|t|$. The hypothesis \eqref{eq:admissible} is a condition on both $t$ and $T$ which ensures that the supremum of $|P|$ on every such segment is at least one twentieth of its supremum on all of $\R^n$. The restricted function $q\mapsto P(qt)$ has domain $\R$; the polynomial $P$ itself has domain $\R^n$.
\begin{proof}
We may assume $P\ne0$. We first average the restriction $q\mapsto P(qt)$ against a Gaussian kernel, show that this averaging almost recovers a fixed multiple of $P$, and then control how far from $u$ a large value can occur.

\emph{Centering the frequencies and estimating complex shifts.}
Center the frequencies by putting
\[
\Psi(\theta)=\e(-\tfrac12\sum_{k=1}^n\theta_k)P(\theta)
=\sum_{\xi\in\{-1/2,1/2\}^n}c_\xi\e(\xi\cdot\theta),\qquad c_\xi=a_{\xi+(1/2,\ldots,1/2)}.
\]
This finite sum defines an entire function on $\C^n$; for $\theta\in\R^n$, $|\Psi(\theta)|=|P(\theta)|$. Here and below $\|\Psi\|_\infty=\sup_{\theta\in\R^n}|\Psi(\theta)|=\|P\|_\infty$. The function $\Psi$ has period two in each coordinate, so the same norm is obtained on $\T_2^n$. For a fixed coordinate and fixed values of all the others, its dependence on that coordinate is $ae^{\pi ix}+be^{-\pi ix}$, with $a,b\in\C$. Since $e^{2\pi ix}$ runs through the unit circle,
\[
\sup_{x\in\R}|ae^{\pi ix}+be^{-\pi ix}|
=\sup_{|\zeta|=1}|a\zeta+b|=|a|+|b|.
\]
For real $y$, the triangle inequality therefore gives
\[
\sup_{x\in\R}|ae^{\pi i(x+iy)}+be^{-\pi i(x+iy)}|
\le |a|e^{-\pi y}+|b|e^{\pi y}
\le e^{\pi|y|}(|a|+|b|).
\]
To keep track of the suprema, fix $y\in\R^n$ and put
$F_j(x)=\Psi(x+i(y_1,\ldots,y_j,0,\ldots,0))$ for $0\le j\le n$ and $x\in\R^n$.
For fixed $x'=(x_k)_{k\ne j}\in\R^{n-1}$, write
$F_{j-1}(x)=a_j(x')e^{\pi ix_j}+b_j(x')e^{-\pi ix_j}$.
The imaginary parts $y_1,\ldots,y_{j-1}$ are fixed within these coefficients; $x'$ consists only of real parts. The one-variable calculation gives
\begin{align*}
\|F_j\|_{\infty,\R^n}
&=\sup_{x'\in\R^{n-1}}
 \left(|a_j(x')|e^{-\pi y_j}+|b_j(x')|e^{\pi y_j}\right)\\
&\le e^{\pi|y_j|}\sup_{x'\in\R^{n-1}}(|a_j(x')|+|b_j(x')|)
=e^{\pi|y_j|}\|F_{j-1}\|_{\infty,\R^n}.
\end{align*}
Iteration from $F_0=\Psi$ to $F_n(x)=\Psi(x+iy)$ proves
\begin{equation}\label{eq:shift}
\sup_{\theta\in\R^n}|\Psi(\theta+iy)|\le e^{\pi\sum_k|y_k|}\|\Psi\|_\infty.
\end{equation}
\emph{The periodized Gaussian and the averaging operator.} Let $\delta=\varepsilon/\sqrt n$, and use normalized Haar measure $m_2$ on $\T_2^n=(\R/2\Z)^n$. The periodized Gaussian density is
\[
G_\delta(\theta)=2^n\sum_{k\in\Z^n}(2\pi\delta^2)^{-n/2}
 e^{-|\theta+2k|^2/(2\delta^2)}.
\]
The Poisson summation formula for Schwartz functions, together with the Gaussian transform, gives the standard absolutely convergent expansion
\begin{equation}\label{eq:heatseries}
G_\delta(\theta)=\sum_{b\in(\frac12\Z)^n}e^{-2\pi^2\delta^2|b|^2}\e(b\cdot\theta).
\end{equation}
These are the usual periodization and convolution formulas; see \cite{RudinGroups}. In particular $G_\delta>0$, $\int G_\delta\dd m_2=1$. For $\xi\in(\frac12\Z)^n$, absolute uniform convergence of \eqref{eq:heatseries} and orthogonality of characters give
\begin{align*}
(G_\delta*\e(\xi\cdot\,))(\theta)
&=\int_{\T_2^n}\sum_{b\in(\frac12\Z)^n}
 e^{-2\pi^2\delta^2|b|^2}\e(b\cdot(\theta-y))\e(\xi\cdot y)\dd m_2(y)\\
&=\sum_b e^{-2\pi^2\delta^2|b|^2}\e(b\cdot\theta)
 \int_{\T_2^n}\e((\xi-b)\cdot y)\dd m_2(y)\\
&=e^{-2\pi^2\delta^2|\xi|^2}\e(\xi\cdot\theta).
\end{align*}
Every frequency of $\Psi$ has $|\xi|^2=n/4$. The multiplier is therefore the same for each term of this finite polynomial, and
\begin{equation}\label{eq:heat}
G_\delta*\Psi=c\Psi,\qquad c=e^{-\pi^2\varepsilon^2/2}=e^{-\pi^2/32}.
\end{equation}
Let $\gamma_{u,T}=N(u,T^2)$ on $\R$. For every $b\in(\frac12\Z)^n$, define
\begin{equation}\label{eq:vb}
v_b=\int\e(-q b\cdot t)\dd\gamma_{u,T}(q)
=\e(-u b\cdot t)e^{-2\pi^2T^2(b\cdot t)^2}.
\end{equation}
The measure $\gamma_{u,T}$ gives Gaussian weight to values of $q$ near $u$, while $G_\delta(\theta-qt)$ weights those points of the line which are near $\theta$ modulo $2\Z^n$. The following operator is intended to reconstruct $c\Psi$ from the values $\Psi(qt)$. Both weights are Gaussian rather than sharp cutoffs; a subsequent moment estimate will locate a large value in a bounded interval.
Define a linear operator $\mathcal T:C(\T_2^n)\to C(\T_2^n)$ by
\[
(\mathcal TF)(\theta)=\int_\R G_\delta(\theta-qt)F(qt)\dd\gamma_{u,T}(q).
\]
In this formula the argument $qt$ is reduced modulo $2\Z^n$ when necessary. In the following sums $\kappa,b\in(\frac12\Z)^n$ and $\xi\in\{-1/2,1/2\}^n$. The kernel series converges absolutely and uniformly and the sum over $\xi$ is finite, so termwise integration is justified. Expanding both factors and then putting $b=\kappa-\xi$ gives
\begin{align*}
\mathcal T\Psi(\theta)
&=\int_\R\left(\sum_\kappa e^{-2\pi^2\delta^2|\kappa|^2}
 \e(\kappa\cdot\theta)\e(-q\kappa\cdot t)\right)
 \left(\sum_\xi c_\xi\e(q\xi\cdot t)\right)\dd\gamma_{u,T}(q)\\
&=\sum_{\kappa,\xi}e^{-2\pi^2\delta^2|\kappa|^2}c_\xi\e(\kappa\cdot\theta)
 \int_\R\e(-q(\kappa-\xi)\cdot t)\dd\gamma_{u,T}(q)\\
&\stackrel{\eqref{eq:vb}}{=}\sum_{\kappa,\xi}e^{-2\pi^2\delta^2|\kappa|^2}c_\xi\e(\kappa\cdot\theta)v_{\kappa-\xi}\\
&=\sum_{b,\xi}e^{-2\pi^2\delta^2|b+\xi|^2}c_\xi
 \e((b+\xi)\cdot\theta)v_b\\
&\stackrel{\eqref{eq:heat}}{=}c\sum_bv_b\e(b\cdot\theta)e^{-2\pi^2\delta^2|b|^2}
 \sum_\xi c_\xi\e(\xi\cdot\theta)e^{-4\pi^2\delta^2b\cdot\xi}.
\end{align*}
We used $|b+\xi|^2=|b|^2+2b\cdot\xi+n/4$. The inner sum is $\Psi(\theta+2\pi i\delta^2b)$, so
\begin{equation}\label{eq:operator}
\mathcal T\Psi(\theta)=c\sum_bv_b\e(b\cdot\theta)e^{-2\pi^2\delta^2|b|^2}\Psi(\theta+2\pi i\delta^2b).
\end{equation}
\emph{A common summable bound.} For real $\theta$, the absolute value of the $b$-summand in \eqref{eq:operator} is bounded as follows:
\begin{align*}
&c|v_b|e^{-2\pi^2\delta^2|b|^2}
       |\Psi(\theta+2\pi i\delta^2b)|\\
&\quad\stackrel{\eqref{eq:shift},\eqref{eq:vb}}{\le}
 c\|\Psi\|_\infty
 \exp\left(-2\pi^2T^2(b\cdot t)^2-2\pi^2\delta^2|b|^2
          +2\pi^2\delta^2\sum_k|b_k|\right)\\
&\quad\le c\|\Psi\|_\infty
 \exp\left(-2\pi^2T^2(b\cdot t)^2-\pi^2\delta^2|b|^2+n\pi^2\delta^2\right).
\end{align*}
The second step uses $2\sum_k|b_k|\le |b|^2+n$. Since $n\delta^2=\varepsilon^2$, weakening the coefficient of the nonpositive $T^2$ term gives the convenient bound
\begin{equation}\label{eq:majorant}
c\|\Psi\|_\infty e^{\pi^2\varepsilon^2}
 e^{-\pi^2\delta^2|b|^2-\pi^2T^2(b\cdot t)^2}.
\end{equation}
This is summable, uniformly in real $\theta$. The $b=0$ term is $c\Psi$, and the other terms give
\begin{equation}\label{eq:operror}
\begin{aligned}
\|\mathcal T\Psi-c\Psi\|_{\infty,\R^n}
&\le c e^{\pi^2\varepsilon^2}D_n(t,T)\|\Psi\|_{\infty,\R^n},\\
\|\mathcal T\Psi\|_{\infty,\R^n}
&\ge c\bigl(1-e^{\pi^2\varepsilon^2}D_n(t,T)\bigr)\|\Psi\|_{\infty,\R^n}
\stackrel{\eqref{eq:admissible}}{\ge}\tfrac c2\|\Psi\|_{\infty,\R^n}.
\end{aligned}
\end{equation}

\emph{Controlling the location of a large value.} Put $H_{u,T}(\theta)=\int G_\delta(\theta-qt)\dd\gamma_{u,T}(q)$ and $D=D_n(t,T)$. The function $H_{u,T}:\R^n\to(0,\infty)$ is two-periodic in each coordinate. Equations \eqref{eq:heatseries} and \eqref{eq:vb} give
\[
H_{u,T}(\theta)=\sum_{b\in(\frac12\Z)^n}
 e^{-2\pi^2\delta^2|b|^2}v_b\e(b\cdot\theta).
\]
The zero-frequency term is one. The sum of the absolute values of all other terms is bounded by $D_n(t,T)$, since its exponential terms have twice the negative exponent used in \eqref{eq:Delta}. Thus
\begin{equation}\label{eq:mass}
|H_{u,T}(\theta)-1|\le D,\qquad |H_{u,\sqrt2T}(\theta)-1|\le D.
\end{equation}
Replacing $T$ by $\sqrt2T$ changes only $v_b$ in this series. By \eqref{eq:vb}, its absolute value $e^{-2\pi^2T^2(b\cdot t)^2}$ can only decrease, so both bounds use the same original $D=D_n(t,T)$.
For $x=(q-u)/T$, the ratio of the densities of $(|q-u|/T)\gamma_{u,T}$ and $\gamma_{u,\sqrt2T}$ is
\[
\frac{\displaystyle\frac{|q-u|}{T}\frac1{\sqrt{2\pi}T}
 e^{-(q-u)^2/(2T^2)}}
{\displaystyle\frac1{\sqrt{2\pi}\sqrt2T}e^{-(q-u)^2/(4T^2)}}
=\sqrt2\frac{|q-u|}{T}e^{-(q-u)^2/(4T^2)}
=\sqrt2|x|e^{-x^2/4}\le2e^{-1/2}<2.
\]
Positivity therefore gives
\begin{equation}\label{eq:moment}
\int\frac{|q-u|}{T}G_\delta(\theta-qt)\dd\gamma_{u,T}(q)
\le2H_{u,\sqrt2T}(\theta)\le2(1+D).
\end{equation}
Finally let $W=\sup_{q\in\R}|\Psi(qt)|/(1+|q-u|/T)$. By its definition and positivity of $G_\delta$,
\begin{align*}
|(\mathcal T\Psi)(\theta)|
&\le\int_\R G_\delta(\theta-qt)|\Psi(qt)|\dd\gamma_{u,T}(q)\\
&\le W\left(H_{u,T}(\theta)+\int_\R\frac{|q-u|}{T}
          G_\delta(\theta-qt)\dd\gamma_{u,T}(q)\right)\\
&\stackrel{\eqref{eq:mass},\eqref{eq:moment}}{\le}3(1+D)W.
\end{align*}
Taking the supremum over $\theta\in\R^n$ gives
\[
\frac c2\|\Psi\|_{\infty,\R^n}
\stackrel{\eqref{eq:operror}}{\le}\|\mathcal T\Psi\|_{\infty,\R^n}
\le3(1+D)W\stackrel{\eqref{eq:admissible}}{\le}6W,
\qquad
\|\Psi\|_{\infty,\R^n}\le\frac{12}{c}W\stackrel{\eqref{eq:heat}}{<}20W.
\]
Since $P\ne0$, also $\|\Psi\|_\infty>0$, and therefore $W>0$. The continuous weighted quotient tends to zero as $|q|\to\infty$, since its numerator is bounded, so its positive maximum is attained. At a maximizing point,
$1+|q-u|/T\le\|\Psi\|_\infty/W<20$ and $|\Psi(qt)|\ge W>\|\Psi\|_\infty/20$. This proves \eqref{eq:sample}.
\end{proof}

\section{Pseudomeasures, uniqueness, and bounded synthesis}\label{sec:PM}
We first work on $\R$, where localization in the small intervals $[x_z,x_z+L_n]$ gives a derivative estimate. Unless a group is explicitly indicated, the pseudomeasure norms in this section are $\PM(\R)$ norms. Supremum norms of transforms are taken over the real frequency variable. We later transfer the conclusion to the circle by a fixed cutoff.
\begin{lemma}\label{lem:derivative}
There is a constant $C_{\rm der}\ge1$ such that every $S\in\PM(\R)$ supported in $[0,L]$, $L>0$, satisfies
\begin{equation}\label{eq:derivative}
\|(\widehat S)'\|_{\infty,\R}\le C_{\rm der}L\|S\|_{\PM}.
\end{equation}
\end{lemma}
\begin{proof}
Fix $h\in C_c^\infty(\R)$ equal to one near $[0,1]$ and put
$g(x)=2\pi ixh(x)$, $C_{\rm der}=\max(1,\|\widehat g\|_1)$. The standard differentiation rule gives
$(\widehat S)'=\calF((2\pi ixh(x/L))S)$. Write $g_L(x)=2\pi ixh(x/L)=Lg(x/L)$. A change of variable gives
\[
\widehat{g_L}(q)=L^2\widehat g(Lq),\qquad
\|\widehat{g_L}\|_{L^1(\R)}=L\|\widehat g\|_{L^1(\R)}.
\]
Consequently, the product--convolution bound \eqref{eq:multiplier} yields
\[
\sup_{q\in\R}|(\widehat S)'(q)|
=\|\widehat{g_LS}\|_{L^\infty(\R)}
\le\|\widehat{g_L}\|_{L^1(\R)}\|S\|_{\PM(\R)}
\le C_{\rm der}L\|S\|_{\PM(\R)}.
\]
\end{proof}
Fix from now on
\begin{equation}\label{eq:tau}
\tau=[10^6(1+C_{\rm der})]^{-1},\qquad T_n=T_n(\tau),\qquad
C_{\rm der}L_nT_n\le\frac{C_{\rm der}\tau}{6}\sqrt{\frac n{n+1}}\le\frac1{6\cdot10^6}.
\end{equation}
The estimate follows from \eqref{eq:tails} and \eqref{eq:T}. Until the end of this section assume $\omega\in\mathcal D_\tau$. Call $n$ admissible if \eqref{eq:admissible} holds with $t=(t_1,\ldots,t_n)$ and $T=T_n$. By the definition of $\mathcal D_\tau$ in Lemma~\ref{lem:Delta}, $\Delta_{2^j}(\omega|_{2^j},\tau)\to0$ as $j\to\infty$. Since the right-hand side of \eqref{eq:admissible} is a fixed positive number, all sufficiently large dyadic indices $n=2^j$ are admissible.

For each fixed $n$ below, write $t=(t_1,\ldots,t_n)\in\R^n$. By the separation in Lemma~\ref{lem:geometry}, for $S\in \PM(\R)$ supported by $E$ and each $n$ we may choose cutoffs $\phi_z\in C_c^\infty(\R)$ equal to one near $E_z$, with pairwise disjoint supports. For $a\in\R$, define the translation operator $\shift_a:\PM(\R)\to\PM(\R)$ by
$\langle\shift_aR,\phi\rangle=\langle R,\phi(\cdot+a)\rangle$, and put $S_z=\shift_{-x_z}(\phi_zS)$. By \eqref{eq:multiplier}, $\|S_z\|_{\PM}\le\|\widehat{\phi_z}\|_1\|S\|_{\PM}$ and $\supp S_z\subset[0,L_n]$. Fix the cutoffs independently of $S$. Replacing each localized piece by an atom at its base point gives the algebraic approximation
\begin{equation}\label{eq:Rdef}
R_nS=\sum_{z\in\{0,1\}^n}\widehat S_z(0)\delta_{x_z}.
\end{equation}
This definition makes sense for every $n$ and every coefficient sequence. The restriction to admissible indices will enter only in the uniform norm estimates. To compare $\widehat{R_nS}$ and $\widehat S$, introduce a real parameter $s$ in the weights:
\begin{equation}\label{eq:Pdef}
P_n:\R^n\times\R\to\C,\quad P_n(\theta,s)=\sum_{z\in\{0,1\}^n}\e(z\cdot\theta)\widehat S_z(s),\quad
\mathcal M_n=\sup_{(\theta,s)\in\R^n\times\R}|P_n(\theta,s)|<\infty.
\end{equation}
Indeed $\mathcal M_n\le\sum_z\|S_z\|_{\PM}$; this preliminary bound need not be uniform in $n$. Since $\widehat{\shift_aR}(q)=\e(aq)\widehat R(q)$ and $x_z=z\cdot t$,
\begin{equation}\label{eq:diagonal}
\widehat S(q)=\calF\left(\sum_z\shift_{x_z}S_z\right)(q)
=\sum_z\e(qx_z)\widehat S_z(q)=P_n(qt,q).
\end{equation}
The other relevant value of $P_n$ is the transform of the atomic approximation. By the finite-sum rule,
\begin{equation}\label{eq:atomicFT}
\widehat{R_nS}(q)
=\sum_z\widehat S_z(0)\widehat{\delta_{x_z}}(q)
=\sum_z\widehat S_z(0)\e(qx_z)
=\sum_z\widehat S_z(0)\e(z\cdot qt)=P_n(qt,0).
\end{equation}
Thus changing the second parameter from $q$ to $0$ is precisely the approximation step. For fixed $\theta\in\R^n$, the pseudomeasure
$S_\theta=\sum_z\e(z\cdot\theta)S_z$ is supported in $[0,L_n]$, and its Fourier transform is $s\mapsto P_n(\theta,s)$. Hence
\[
\|S_\theta\|_{\PM(\R)}=\sup_{s\in\R}|P_n(\theta,s)|\le\mathcal M_n.
\]
Applying Lemma~\ref{lem:derivative} to this pseudomeasure gives
\begin{equation}\label{eq:Pderivative}
\begin{aligned}
|\partial_sP_n(\theta,s)|
&=|(\widehat{S_\theta})'(s)|\le\|(\widehat{S_\theta})'\|_{\infty,\R}\\
&\stackrel{\text{Lemma }\ref{lem:derivative}}{\le}
C_{\rm der}L_n\|S_\theta\|_{\PM(\R)}
\le C_{\rm der}L_n\mathcal M_n.
\end{aligned}
\end{equation}
Integration with $\theta=qt$ fixed gives the frequently used comparison
\begin{equation}\label{eq:compare}
|\widehat S(q)-P_n(qt,x)|
\stackrel{\eqref{eq:diagonal}}{=}\left|\int_x^q\partial_sP_n(qt,s)\dd s\right|
\le C_{\rm der}L_n\mathcal M_n|q-x|\qquad(q,x\in\R).
\end{equation}
\begin{lemma}\label{lem:M}
For every admissible $n$, $\mathcal M_n\le80\|S\|_{\PM(\R)}$.
\end{lemma}
\begin{proof}
If $\mathcal M_n=0$, the conclusion is immediate. Otherwise choose $(\theta_0,s_0)\in\R^n\times\R$ with $|P_n(\theta_0,s_0)|\ge\mathcal M_n/2$. The polynomial $P_n(\cdot,s_0):\R^n\to\C$ has supremum norm at least $\mathcal M_n/2$. Lemma~\ref{lem:sampling}, with center $s_0$, therefore gives $|q-s_0|\le20T_n$ and $|P_n(qt,s_0)|\ge\mathcal M_n/40$. Equations \eqref{eq:compare} and \eqref{eq:tau} imply
\[
\begin{aligned}
\|S\|_{\PM(\R)}&\ge|\widehat S(q)|\\
&\ge |P_n(qt,s_0)|-|\widehat S(q)-P_n(qt,s_0)|\\
&\stackrel{\eqref{eq:compare}}{\ge}\frac{\mathcal M_n}{40}
 -C_{\rm der}L_n\mathcal M_n|q-s_0|\\
&\stackrel{\eqref{eq:tau}}{\ge}
\left(\frac1{40}-\frac{20}{6\cdot10^6}\right)\mathcal M_n
\ge\frac{\mathcal M_n}{80}.
\end{aligned}
\]
\end{proof}
\begin{proposition}\label{prop:tail}
For every coefficient sequence in $\mathcal D_\tau$, inequality \eqref{eq:tailgoal} holds for all pseudomeasures supported by $E$.
\end{proposition}
\begin{proof}
Let $M=\|S\|_{\PM}>0$ and choose $q_0$ with $|\widehat S(q_0)|>3M/4$. Fix $R>0$ and put $u_n=R+21T_n$. Along sufficiently large admissible dyadic indices, we compare directly with $\widehat S(q_0)=P_n(q_0t,q_0)$:
\[
\begin{aligned}
|P_n(q_0t,u_n)-\widehat S(q_0)|
&\stackrel{\eqref{eq:compare}}{\le}C_{\rm der}L_n\mathcal M_n|u_n-q_0|\\
&\stackrel{\text{Lemma }\ref{lem:M}}{\le}80C_{\rm der}L_nM(R+|q_0|+21T_n)\\
&\stackrel{\eqref{eq:tau}}{\le}80C_{\rm der}L_nM(R+|q_0|)
 +\frac{1680}{6\cdot10^6}M\le\frac M8.
\end{aligned}
\]
The last inequality holds for sufficiently large $n$, since $R,q_0$ are fixed and $L_n\to0$ by \eqref{eq:tails}. Therefore
\[
\begin{aligned}
\|P_n(\cdot,u_n)\|_\infty
&\ge |P_n(q_0t,u_n)|\\
&\ge |\widehat S(q_0)|-|P_n(q_0t,u_n)-\widehat S(q_0)|\\
&>\frac{3M}{4}-\frac M8=\frac{5M}{8}\ge\frac M2.
\end{aligned}
\]
Lemma~\ref{lem:sampling} gives $q\ge u_n-20T_n>R$ with $|P_n(qt,u_n)|\ge M/40$. At this point the two arguments of $P_n$ still differ. Equations \eqref{eq:compare} and \eqref{eq:tau}, together with Lemma~\ref{lem:M}, control the error of replacing $u_n$ by $q$:
\[
\begin{aligned}
|\widehat S(q)-P_n(qt,u_n)|
&\stackrel{\eqref{eq:compare}}{\le} C_{\rm der}L_n\mathcal M_n|q-u_n|\\
&\stackrel{\text{Lemma }\ref{lem:M}}{\le} C_{\rm der}L_n(80M)(20T_n)
\stackrel{\eqref{eq:tau}}{\le}\frac{1600M}{6\cdot10^6}<\frac M{80}.
\end{aligned}
\]
It follows that $|\widehat S(q)|\ge M/40-M/80=M/80$ for some $q>R$. Since $R>0$ was arbitrary, $\limsup_{q\to+\infty}|\widehat S(q)|\ge M/80$, as required.
\end{proof}

We retain a short localization argument connecting the line and the circle, to make explicit the norm and decay needed here.
\begin{lemma}\label{lem:lift}
Every $S_\T\in\PM(\T)$ supported in $E\subset[0,2/15]$ has a lift $S_\R\in\PM(\R)$ supported in the same real set, with
\[
\widehat S_\R(k)=\widehat S_\T(k)\quad(k\in\Z),\qquad
\|S_\R\|_{\PM(\R)}\le C_\chi\|S_\T\|_{\PM(\T)}.
\]
The constant $C_\chi$ is fixed, independent of $E$ and $S_\T$. If $S_\T\in\PF(\T)$, then $S_\R\in\PF(\R)$.
\end{lemma}
\begin{proof}
Choose $\chi\in C_c^\infty((-1/4,1/2))$ equal to one near $[0,2/15]$. For $\Pi:C_c^\infty(\R)\to C^\infty(\T)$, $\Pi f(x)=\sum_{m\in\Z}f(x+m)$, define
$\langle S_\R,\phi\rangle=\langle S_\T,\Pi(\chi\phi)\rangle$.
The supporting interval has length less than one, so this distribution is supported in the real copy of $E$. For fixed $s\in\R$, put $F_s=\Pi(\chi\e(s\cdot))$. The coefficient of $\e(kx)$ in its Fourier series is
\[
\begin{aligned}
\int_0^1F_s(x)\e(-kx)\dd x
&=\sum_{\ell\in\Z}\int_0^1
 \chi(x+\ell)\e(s(x+\ell))\e(-kx)\dd x\\
&=\sum_{\ell\in\Z}\int_\ell^{\ell+1}
 \chi(y)\e((s-k)y)\dd y\\
&=\int_\R\chi(y)\e((s-k)y)\dd y
 =\widehat\chi(s-k).
\end{aligned}
\]
Here we used $y=x+\ell$ and $\e(k\ell)=1$ for $k,\ell\in\Z$; compact support of $\chi$ justifies the interchange of sum and integral. Rapid decay of $\widehat\chi$ implies that
\[
F_s=\sum_{k\in\Z}\widehat\chi(s-k)\e(k\cdot)
\quad\text{with convergence in }A(\T).
\]
Since $\widehat S_\R(s)=\langle S_\T,F_s\rangle$, applying the continuous functional $S_\T$ term by term yields the first identity below. The uniform bound stated alongside it will be verified next:

\begin{equation}\label{eq:lift}
\widehat S_\R(s)=\sum_{k\in\Z}\widehat\chi(s-k)\widehat S_\T(k),\qquad
C_\chi=\sup_{s\in\R}\sum_{k\in\Z}|\widehat\chi(s-k)|<\infty.
\end{equation}
Since $\chi\in C_c^\infty(\R)$, integration by parts twice gives $|\widehat\chi(v)|\le C(1+|v|)^{-2}$. Write $s=a+r$ with $a\in\Z$ and $r\in[-1/2,1/2]$, and shift the summation index. Then
\[
\sum_{k\in\Z}|\widehat\chi(s-k)|
=\sum_{k\in\Z}|\widehat\chi(r-k)|
\le C\sum_{k\in\Z}(1+|r-k|)^{-2}
\le4C\sum_{k\in\Z}(1+|k|)^{-2}<\infty.
\]
The uniform bound justifies absolute convergence of \eqref{eq:lift} and gives
\[
|\widehat S_\R(s)|
\le\sum_{k\in\Z}|\widehat\chi(s-k)|\,|\widehat S_\T(k)|
\le C_\chi\|S_\T\|_{\PM(\T)},\qquad
\|S_\R\|_{\PM(\R)}\le C_\chi\|S_\T\|_{\PM(\T)}.
\]
For $m\in\Z$, pairings of the compactly supported distribution $S_\R$ with $\e(m\cdot)$ mean pairings after multiplication by a smooth compactly supported cutoff equal to one near $\supp S_\R$. We may choose that cutoff equal to one on $\supp\chi$, so that the definition of $S_\R$ gives
\begin{align*}
\widehat S_\R(m)
&=\langle S_\R,\e(m\cdot)\rangle
=\langle S_\T,\Pi(\chi\e(m\cdot))\rangle\\
&=\langle S_\T,\e(m\cdot)\Pi\chi\rangle
=\langle S_\T,\e(m\cdot)\rangle=\widehat S_\T(m).
\end{align*}
Here $\Pi(\chi\e(m\cdot))=\e(m\cdot)\Pi\chi$ because $m\in\Z$, and the next equality uses $\Pi\chi=1$ near $\supp S_\T$. If $\widehat S_\T(k)\to0$, split \eqref{eq:lift} at $|k|\le K$. The resulting bound is
\[
|\widehat S_\R(s)|\le
\sum_{|k|\le K}|\widehat\chi(s-k)|\,|\widehat S_\T(k)|
+C_\chi\sup_{|k|>K}|\widehat S_\T(k)|.
\]
For fixed $K$, the finite sum tends to zero as $|s|\to\infty$. Taking the upper limit and then letting $K\to\infty$ proves $\widehat S_\R(s)\to0$.
\end{proof}
\begin{corollary}\label{cor:U}
For $\omega\in\mathcal D_\tau$, $E\in U$.
\end{corollary}
\begin{proof}
Suppose that $0\ne S_\T\in\PF(\T)$ is supported by $E$, and let $S_\R\in\PF(\R)$ be its lift from Lemma~\ref{lem:lift}. By uniqueness of Fourier coefficients on the circle, some $m\in\Z$ satisfies
\[
\widehat S_\R(m)=\widehat S_\T(m)\ne0.
\]
Since $\widehat S_\R$ is continuous, its absolute value is bounded below by a positive constant on a neighborhood of $m$. Hence its essential supremum is positive, so $\|S_\R\|_{\PM(\R)}>0$. On the other hand, $S_\R\in\PF(\R)$ implies $\widehat S_\R(q)\to0$ as $q\to+\infty$, contradicting Proposition~\ref{prop:tail}. Apply \eqref{eq:U}.
\end{proof}

For $G=\R$ or $G=\T$, write $\PM_G(E)=\{S\in \PM(G):\supp S\subset E\}$, with the inherited pseudomeasure norm. On the circle this is the space denoted $\PM(E)$ in \eqref{eq:U}. A linear operator has finite rank when its range is finite-dimensional.
\begin{proposition}\label{prop:synthesis}
For $\omega\in\mathcal D_\tau$, along the sufficiently large dyadic indices there are linear finite-rank operators $R_n:\PM_\R(E)\to \PM_\R(E)$, with atomic values supported by $E$, such that
\begin{gather}
\|R_nS\|_{\PM(\R)}\le80\|S\|_{\PM(\R)},\qquad R_nS\xrightarrow{w^*}S,\label{eq:synbound}\\
|\widehat{R_nS}(q)-\widehat S(q)|\le80C_{\rm der}L_n|q|\|S\|_{\PM(\R)}.\label{eq:synerror}
\end{gather}
Here weak-star convergence is in $\PM(\R)=A(\R)^*$, that is, convergence of $\langle R_nS,f\rangle$ for every $f\in A(\R)$. Analogous operators on $\PM(E)\subset\PM(\T)$ have norm at most $80C_\chi$ and converge weak-star in $A(\T)^*$. The set $E$ is a set of spectral synthesis on both groups.
\end{proposition}
\begin{proof}
Use the operators $R_n$ defined in \eqref{eq:Rdef}, with the cutoffs already fixed independently of $S$. Its range lies in $\operatorname{span}\{\delta_{x_z}:z\in\{0,1\}^n\}$, so $\rank R_n\le2^n$. Since $x_z\in E$ and the calculation \eqref{eq:atomicFT} gives $\widehat{R_nS}(q)=P_n(qt,0)$, Lemma~\ref{lem:M} and \eqref{eq:compare} give the estimates. Pointwise convergence of the transforms and their uniform bound imply weak-star convergence against $L^1(\R)$ by dominated convergence. On the circle compose with the fixed lift of Lemma~\ref{lem:lift} and project the atomic measures back; the norm bound is $80C_\chi$, and convergence follows at every integer coefficient.
For either group let $I(E)=\{f\in A:f|_E=0\}$ and let $J(E)$ be the norm closure of the functions in $A$ vanishing on a neighborhood of $E$. The standard support duality \cite{RudinGroups,KL} states precisely that
\[
J(E)^\perp=\PM_G(E).
\]
Our goal is to prove $\langle S,f\rangle=0$ for every $S\in\PM_G(E)$ and every $f\in I(E)$. Here $J(E)^\perp$ denotes the annihilator in the Banach dual $A^*$, rather than the measure polar used in the introduction. If $f\in I(E)$ and $S\in J(E)^\perp$, then
\[
\langle S,f\rangle
=\lim_n\langle R_nS,f\rangle
=\lim_n\sum_z\widehat S_z(0)f(x_z)=0,
\]
because every $x_z$ belongs to $E$. (On the circle use the corresponding projected operators.) Thus $J(E)^\perp\subset I(E)^\perp$. Since $J(E)\subset I(E)$, the opposite inclusion is automatic, and Hahn--Banach yields $I(E)=J(E)$. This equality is the definition of spectral synthesis.
\end{proof}
The bounds here are pseudomeasure norms, not total variation bounds for the atomic measures. Neither norm convergence nor positivity of the approximants for a general pseudomeasure is asserted.

For later use, consider the same algebraic prescription on a finite positive measure $\nu$ supported by $E$, without assuming $\omega\in\mathcal D_\tau$. Then
\begin{equation}\label{eq:positiveRn}
\widehat{\nu_z}(0)=\nu(E_z),\qquad
R_n\nu=\sum_{z\in\{0,1\}^n}\nu(E_z)\delta_{x_z}
\quad(n\ge1).
\end{equation}
Its mass is $\nu(E)$, and for every continuous test function $f$,
\[
\left|\int f\dd(R_n\nu)-\int f\dd\nu\right|
\le\nu(E)\sup\{|f(x)-f(y)|:x,y\in E,\ |x-y|\le L_n\}\longrightarrow0.
\]
Thus these positive measures converge weak-star for all $n\to\infty$, independently of the uniform pseudomeasure bounds along admissible indices. In Section~\ref{sec:full} we apply a sequence of maps to these approximating measures. Since the maps themselves vary along the sequence, weak-star convergence alone does not control the resulting image measures; a separate quantitative estimate is needed.

\section{Supports and quotients of tori}\label{sec:groups}
We collect the standard compact-group facts in the precise form needed below; the duality input is in \cite{RudinGroups}. The topology of $\T^d$ is the product topology of $d$ circles. For $\Gamma\le\Z^d$, write
\[
\Gamma^\perp=\{x\in\T^d:\e(\alpha\cdot x)=1\text{ for every }\alpha\in\Gamma\}.
\]
For a subgroup $G\le\T^d$, conversely, set
\[
G^\perp=\{\alpha\in\Z^d:\e(\alpha\cdot x)=1\text{ for every }x\in G\}.
\] Closed subgroups satisfy $G=(G^\perp)^\perp$. Every subgroup of $\Z^d$ admits a finite $\Z$-basis, which is also linearly independent over $\R$ (rational matrices have the same rank over $\Q$ and $\R$).

\begin{lemma}\label{lem:groups}
Let $A\in\Z^{r\times d}$ have rank $r\ge1$, and let $\sigma$ be a finite positive measure on $\T^d$, with $F=\supp\sigma$.
Then $A:\T^d\to\T^r$ is onto, $\supp(A_*\sigma)=A(F)$, and
\begin{equation}\label{eq:fibres}
F+\ker A=F\quad\Longrightarrow\quad F=A^{-1}(\supp(A_*\sigma)).
\end{equation}
If $\widehat\sigma$ vanishes outside a subgroup $\Gamma\le\Z^d$, then $\sigma$ is translation invariant under $\Gamma^\perp$. If the rows of $A$ are a basis of $\Gamma$, $\Gamma^\perp=\ker A$.
For nonzero positive measures $\kappa,\gamma$,
\begin{equation}\label{eq:convsupp}
\supp(\kappa*\gamma)=\supp\kappa+\supp\gamma.
\end{equation}
\end{lemma}
\begin{proof}
Real surjectivity gives torus surjectivity, and integrality ensures the map is well-defined. The continuous image of a compact support is the support of the image measure: inverse images of neighborhoods of image points have positive measure. Since $A^{-1}(A(F))=F+\ker A$, this proves \eqref{eq:fibres}. For $h\in\Gamma^\perp$, translating the measure $\sigma$ by $h$, that is, taking its image under $x\mapsto x+h$, multiplies its Fourier coefficient at $\alpha$ by $\e(\alpha\cdot h)$. The assumed vanishing makes all coefficients unchanged for $h\in\Gamma^\perp$; uniqueness of Fourier coefficients proves invariance. The basis description gives the kernel identity. Finally, let $a(x,y)=x+y$. Convolution is the image $\kappa*\gamma=a_*(\kappa\otimes\gamma)$, and hence
\[
\supp(\kappa*\gamma)=a(\supp(\kappa\otimes\gamma))=a(\supp\kappa\times\supp\gamma)=\supp\kappa+\supp\gamma.
\]
\end{proof}
The kernels here are on tori. Even when $r=d$, $\ker A$ need not be trivial: multiplication by two on $\T$ has kernel $\{0,1/2\}$.

\begin{lemma}\label{lem:gaussgroup}
Let $\Sigma\ne0$ be a symmetric positive semidefinite $d\times d$ matrix, and let $\gamma=(\rho_d)_*N(0,\Sigma/4)$, with $\rho_d$ the canonical projection onto $\T^d$ defined in Section~\ref{sec:prob}. Then
\[
G=\supp\gamma=\overline{\rho_d(\operatorname{range}\Sigma)}
\]
is a nontrivial connected closed subgroup, and
\begin{equation}\label{eq:annihilator}
\Lambda=G^\perp=\{\alpha\in\Z^d:\Sigma\alpha=0\},\qquad \rank\Lambda<d.
\end{equation}
If $r=\rank\Lambda>0$ and $C_\Lambda$ has as rows a basis of $\Lambda$, then $G=\ker C_\Lambda$. If $r=0$, then $G=\T^d$.
\end{lemma}
\begin{proof}
The Gaussian has full relative support in $V=\operatorname{range}\Sigma$. Its projected support is the closure of the image of this connected subgroup. Since $V\ne\{0\}$, that image is nontrivial. A character is one on the image of $V$ exactly when its frequency is orthogonal to $V$, equivalently $\Sigma\alpha=0$. Thus $\rank\Lambda\le\dim\ker_\R\Sigma<d$. If $r>0$, then
\[
\ker C_\Lambda
\overset{\text{Lemma~\ref{lem:groups}}}{=}\Lambda^\perp
=(G^\perp)^\perp=G,
\]
where the last equality follows from the closedness of $G$ and the double-annihilator identity recalled above. If $r=0$, then $\Lambda=\{0\}$, so the same identity gives $G=\Lambda^\perp=\T^d$.
\end{proof}
The closure is taken in $\T^d$ and is essential: a real subspace may have a nonclosed image in the torus.

\section{Full support of limiting images}\label{sec:full}
Throughout this section fix $\omega\in\mathcal R$, where $\mathcal R$ is the event of Lemma~\ref{lem:energy}. Thus \eqref{eq:robust} holds for all $n\ge n_*(\omega)$. The coefficients $t_k=\omega_k/Q_k$ are now fixed; only the digits on $(\mathcal Z,\PP_{\rm dig})$ remain random.
\begin{theorem}\label{thm:full}
For every $d\ge1$, every quasi-independent sequence $(a^j)$ in $\Z^d$, and every nonzero finite positive measure $\nu\ll\mu$, all weak-star cluster points of $(m_{a^j})_*\nu$ have support $\T^d$.
\end{theorem}
Fix $d$, $(a^j)$, and $\nu$ as in the theorem, and put $m=\nu(E)>0$. All image measures have mass $m$ and live on the compact metrizable space $\T^d$. The set of positive measures of mass $m$ is weak-star compact and metrizable, so cluster points exist and each is the limit of a subsequence. Choose any such point $\sigma$, pass to that subsequence, and relabel so that
\begin{equation}\label{eq:fixedlimit}
(m_{a^j})_*\nu\xrightarrow{w^*}\sigma,\qquad \sigma(\T^d)=m.
\end{equation}
The following three lemmas use this fixed data; every further subsequence retains the same limit. They do not assume the full-support conclusion in any dimension. The final proof will use induction only after these preparatory lemmas have been established.

The proof first separates the finite prefix carrying the density from the remaining independent digits. The Fourier coefficients of those remaining digits define a subgroup $\Gamma\le\Z^d$ of bounded energy. If $\rank\Gamma<d$, translation invariance reduces the problem immediately to a lower-dimensional torus. If $\rank\Gamma=d$, the deletion estimate produces a nonzero Gaussian factor; its support then supplies a lower-dimensional quotient. The next three lemmas implement these steps.

\begin{lemma}[Removing the tail and approximating the density]\label{lem:prefix}
Put $B_j=\max_\ell|a^j_\ell|$. Discarding finitely many initial indices, choose $n_j$ by $Q_{n_j}\le B_j<Q_{n_j+1}$ and set
\begin{equation}\label{eq:indices}
h_j=\lfloor\sqrt{n_j}\rfloor,\qquad N_j=n_j+2,\qquad F_j=\{h_j+1,\ldots,N_j\}.
\end{equation}
Define measures on $\T^d$ by
\begin{equation}\label{eq:betaeta}
\begin{aligned}
\beta_j&=\sum_{z\in\{0,1\}^{h_j}}\nu(E_z)\delta_{\rho_d(a^jx_z)},\\
\eta_j&=2^{-|F_j|}\sum_{\zeta\in\{0,1\}^{F_j}}
 \delta_{\rho_d(a^j\sum_{k\in F_j}t_k\zeta_k)}.
\end{aligned}
\end{equation}
Then $\beta_j$ has mass $m$, $\eta_j$ is a probability measure, and
\begin{equation}\label{eq:prefixlimit}
\beta_j*\eta_j\xrightarrow{w^*}\sigma,\qquad
\widehat\beta_j(\alpha)\widehat\eta_j(\alpha)\longrightarrow\widehat\sigma(\alpha)
\quad(\alpha\in\Z^d).
\end{equation}
\end{lemma}
Here $\{0,1\}^{F_j}$ denotes the binary families $\zeta=(\zeta_k)_{k\in F_j}$ indexed by $F_j$. The base points $x_z$ and pieces $E_z$ are those of Lemma~\ref{lem:geometry}. By \eqref{eq:positiveRn}, $\beta_j=(m_{a^j})_*(R_{h_j}\nu)$. This uses the algebraic prescription for arbitrary $h_j$, not the uniform bounds restricted to admissible indices. For the canonical measure $\mu$, the finite binary product gives the exact identity
\begin{equation}\label{eq:canonicalsplit}
(m_{a^j})_*(R_{N_j}\mu)=(m_{a^j})_*(R_{h_j}\mu)*\eta_j.
\end{equation}
For general $\nu\ll\mu$, averaging its density over prefix cylinders supplies the corresponding approximation.
\begin{proof}
Quasi-independence implies $B_j\to\infty$, so $n_j\to\infty$ and $h_j\to\infty$. Recall the homeomorphism $\pi:\mathcal Z\to E$, $\pi(z)=\sum_kt_kz_k$, with $\mu=\pi_*\PP_{\rm dig}$. Choose a nonnegative density and put $f=(\dd\nu/\dd\mu)\circ\pi\in L^1(\PP_{\rm dig})$.
For $h\ge1$, define cylinder averaging $C_h:L^1(\PP_{\rm dig})\to L^1(\PP_{\rm dig})$ by
\[
(C_hg)(z)=2^h\int_{[z|_h]}g\dd\PP_{\rm dig}.
\]
It preserves integrals, fixes functions depending only on the first $h$ digits, and is a contraction:
\[
\|C_hg\|_1
=\sum_{\zeta\in\{0,1\}^h}\left|\int_{[\zeta]}g\dd\PP_{\rm dig}\right|
\le\sum_{\zeta}\int_{[\zeta]}|g|\dd\PP_{\rm dig}=\|g\|_1.
\]
Set
\begin{equation}\label{eq:densityprefix}
f_j=C_{h_j}f,\qquad
f_j(z)=2^{h_j}\int_{[z|_{h_j}]}f\dd\PP_{\rm dig}
=2^{h_j}\nu(E_{z|_{h_j}}).
\end{equation}
Thus $f_j\ge0$ and $\int f_j\dd\PP_{\rm dig}=m$. Functions depending on finitely many digits are dense in $L^1(\PP_{\rm dig})$. Given $a>0$, choose such a $g$ with $\|f-g\|_1<a$. For all large $j$, $C_{h_j}g=g$, and
\[
\|f_j-f\|_1\le\|C_{h_j}(f-g)\|_1+\|g-f\|_1\le2\|f-g\|_1<2a.
\]
All norms in this paragraph are in $L^1(\PP_{\rm dig})$.

Define maps $X_j,X_j^{\le N_j}:\mathcal Z\to\T^d$ by
\[
X_j(z)=\rho_d(a^j\pi(z)),\qquad
X_j^{\le N_j}(z)=\rho_d\left(a^j\sum_{k\le N_j}t_kz_k\right).
\]
The finite product of the prefix and the block $F_j$ gives
\begin{align*}
(m_{a^j})_*\nu&=(X_j)_*(f\PP_{\rm dig}),\\
(m_{a^j})_*(R_{N_j}\nu)&=(X_j^{\le N_j})_*(f\PP_{\rm dig}),\\
\beta_j*\eta_j&=(X_j^{\le N_j})_*(f_j\PP_{\rm dig}).
\end{align*}
For the last identity, the weight of a prefix $z$ and a block $\zeta$ is
$2^{h_j}\nu(E_z)2^{-N_j}=\nu(E_z)2^{-|F_j|}$, exactly the weight in the convolution of \eqref{eq:betaeta}. Taking $f=f_j=1$ also proves \eqref{eq:canonicalsplit}.

\emph{Change of density.} Images contract total variation, so
\begin{equation}\label{eq:densityTV}
\left\|(m_{a^j})_*(R_{N_j}\nu)-\beta_j*\eta_j\right\|_{\rm TV}
\le\|f-f_j\|_{L^1(\PP_{\rm dig})}\longrightarrow0.
\end{equation}
\emph{Removal of the tail.} The growth scale and \eqref{eq:tails} give
\begin{equation}\label{eq:imagetrail}
B_jL_{N_j}\le\frac{Q_{n_j+1}}{6Q_{n_j+2}\sqrt{n_j+3}}
=\frac1{6M_{n_j+2}\sqrt{n_j+3}}\longrightarrow0.
\end{equation}
Use the quotient Euclidean distance $d_{\T^d}(\rho_d(x),\rho_d(y))=\inf_{k\in\Z^d}|x-y-k|$. Then
\[
d_{\T^d}(X_j(z),X_j^{\le N_j}(z))
\le |a^j|\sum_{k>N_j}t_kz_k\le\sqrt d B_jL_{N_j}.
\]
For $\varphi\in C(\T^d)$, put
$\operatorname{osc}_\varphi(r)=\sup\{|\varphi(x)-\varphi(y)|:d_{\T^d}(x,y)\le r\}$.
It follows that
\[
\left|\int\varphi\dd(m_{a^j})_*\nu-\int\varphi\dd(m_{a^j})_*(R_{N_j}\nu)\right|
\le m\operatorname{osc}_\varphi(\sqrt d B_jL_{N_j})\longrightarrow0.
\]
Combining this weak-star error with \eqref{eq:densityTV} yields the quantitative estimate
\begin{equation}\label{eq:blockapprox}
\begin{aligned}
&\left|\int\varphi\dd(m_{a^j})_*\nu-\int\varphi\dd(\beta_j*\eta_j)\right|\\
&\qquad\le\|\varphi\|_{\infty,\T^d}\|f-f_j\|_{L^1(\PP_{\rm dig})}
 +m\operatorname{osc}_\varphi(\sqrt d B_jL_{N_j})\longrightarrow0.
\end{aligned}
\end{equation}
Here the modulus tends to zero by uniform continuity on the compact torus. Together with \eqref{eq:fixedlimit}, this proves \eqref{eq:prefixlimit}, including its Fourier-coefficient assertion.
\end{proof}

\begin{lemma}[The finite-energy subgroup]\label{lem:finiteenergy}
Use the measures and index sets of Lemma~\ref{lem:prefix}. For $\alpha\in\Z^d$ define
\[
E_j(\alpha)=\sum_{k\in F_j}\tn{(\alpha\cdot a^j)t_k}^2.
\]
There is a strictly increasing sequence $(j_\ell)$ such that every $E_{j_\ell}(\alpha)$ converges in $[0,\infty]$. For this subsequence,
\[
\Gamma=\{\alpha\in\Z^d:\lim_{\ell\to\infty}E_{j_\ell}(\alpha)<\infty\}
\]
is a subgroup of $\Z^d$, and $\widehat\sigma(\alpha)=0$ for $\alpha\notin\Gamma$.
\end{lemma}
\begin{proof}
Compactness of $[0,\infty]$ and a diagonal selection over the countable set $\Z^d$ give $(j_\ell)$. The triangle inequality on the circle and Minkowski's inequality give
$\sqrt{E_j(\alpha+\beta)}\le\sqrt{E_j(\alpha)}+\sqrt{E_j(\beta)}$.
Together with $E_j(0)=0$ and $E_j(-\alpha)=E_j(\alpha)$, this proves that $\Gamma$ is a subgroup. Directly from the finite Dirac sum \eqref{eq:betaeta},
\begin{align*}
\widehat\eta_j(\alpha)
&=2^{-|F_j|}\sum_{\zeta\in\{0,1\}^{F_j}}
 \e\left((\alpha\cdot a^j)\sum_{k\in F_j}t_k\zeta_k\right)\\
&=2^{-|F_j|}\sum_{\zeta}\prod_{k\in F_j}\e((\alpha\cdot a^j)t_k\zeta_k)\\
&=\prod_{k\in F_j}\left(\frac12\sum_{r=0}^1\e((\alpha\cdot a^j)t_kr)\right)
=\prod_{k\in F_j}\frac{1+\e((\alpha\cdot a^j)t_k)}2.
\end{align*}
Since $1+\e(x)=2e^{\pi ix}\cos(\pi x)$, taking absolute values gives
\begin{equation}\label{eq:cosine}
|\widehat\eta_j(\alpha)|=
\prod_{k\in F_j}|\cos(\pi(\alpha\cdot a^j)t_k)|\le e^{-2E_j(\alpha)}.
\end{equation}
The scalar inequality $|\cos(\pi x)|\le e^{-2\tn{x}^2}$ follows from concavity of $\log\cos(\pi x)+2x^2$ on $(-1/2,1/2)$ and periodicity. Since $\beta_j$ is a positive measure of mass $m$ by Lemma~\ref{lem:prefix}, we have $|\widehat\beta_j(\alpha)|\le\beta_j(\T^d)=m$ for every $\alpha\in\Z^d$. If $\alpha\notin\Gamma$, the selected limit is infinite. Hence
\[
|\widehat\beta_{j_\ell}(\alpha)\widehat\eta_{j_\ell}(\alpha)|
\le m e^{-2E_{j_\ell}(\alpha)}\longrightarrow0.
\]
Lemma~\ref{lem:prefix}, specifically \eqref{eq:prefixlimit}, identifies this limit with $\widehat\sigma(\alpha)$.
\end{proof}
From now on in the proof, pass to this subsequence and relabel its indices by $j$. In particular, $\sup_jE_j(\alpha)<\infty$ for $\alpha\in\Gamma$, while the weak-star limit remains $\sigma$.

\begin{lemma}[A nonzero Gaussian factor]\label{lem:factor}
In Lemma~\ref{lem:finiteenergy}, suppose $\rank\Gamma=d$ and let $A\in\Z^{d\times d}$ have as rows a basis of $\Gamma$. There exist a nonzero finite positive Borel measure $\kappa$ on $\T^d$ and a nonzero symmetric positive semidefinite matrix $\Sigma$ such that
\begin{equation}\label{eq:factor}
A_*\sigma=\kappa*\gamma,
\end{equation}
where $\gamma=(\rho_d)_*N(0,\Sigma/4)$.
\end{lemma}
\begin{proof}
\emph{Keeping a nonzero amount of energy in small increments.} Put $b^j=Aa^j$. For $k\in F_j$, let
$v_{j,k}=(\pr{b^j_1t_k},\ldots,\pr{b^j_dt_k})$, where $\pr{x}$ is the unique representative in $[-1/2,1/2)$. Write the rows of $A$ as $\alpha_1,\ldots,\alpha_d\in\Gamma$. Then
\begin{equation}\label{eq:energybound}
\begin{aligned}
\sum_{k\in F_j}|v_{j,k}|^2
&=\sum_{k\in F_j}\sum_{\ell=1}^d\tn{b^j_\ell t_k}^2
=\sum_{\ell=1}^d\sum_{k\in F_j}\tn{(\alpha_\ell\cdot a^j)t_k}^2\\
&=\sum_{\ell=1}^d E_j(\alpha_\ell)
\le\sum_{\ell=1}^d\sup_jE_j(\alpha_\ell)=:C<\infty.
\end{aligned}
\end{equation}
The sequence $b^j$ is quasi-independent, since $\alpha\cdot b^j=(A^{\mathsf T}\alpha)\cdot a^j$ and $A^{\mathsf T}\alpha\in\Z^d\setminus\{0\}$ whenever $\alpha\in\Z^d\setminus\{0\}$.
Recall from Lemma~\ref{lem:prefix} that
\[
B_j=\max_{1\le\ell\le d}|a^j_\ell|,\qquad Q_{n_j}\le B_j<Q_{n_j+1}.
\]
Choose $\ell_j$ maximizing $|b^j_\ell|$, put $k_j=b^j_{\ell_j}\in\Z$, and choose $n'_j$ with $Q_{n'_j}\le|k_j|<Q_{n'_j+1}$. There are fixed $c,C_A>0$ with $cB_j\le|k_j|\le C_AB_j$; one may take $c=\|A^{-1}\|_{\infty\to\infty}^{-1}$ and $C_A=\|A\|_{\infty\to\infty}$. Here $\|A\|_{\infty\to\infty}=\max_\ell\sum_k|A_{\ell k}|$ is the matrix norm induced by the maximum norm on vectors. For large $j$, $cM_{n_j}\ge1$ and $M_{n_j+2}>C_A$, whence
\begin{align*}
|k_j|&\ge cB_j\ge cQ_{n_j}=cM_{n_j}Q_{n_j-1}\ge Q_{n_j-1},\\
|k_j|&\le C_AB_j<C_AQ_{n_j+1}<M_{n_j+2}Q_{n_j+1}=Q_{n_j+2}.
\end{align*}
Thus $n_j-1\le n'_j\le n_j+1<N_j$, and $n'_j\to\infty$.

Let $\delta_j=(n'_j)^{-1/8}$ and $I_j=\{k\in F_j:|v_{j,k}|\le\delta_j\}$. By \eqref{eq:energybound}, at most $C\delta_j^{-2}$ indices in $F_j$ are discarded. Among $1,\ldots,n'_j$, the total discarded number is at most
$h_j+C\delta_j^{-2}=O((n'_j)^{1/2})+O((n'_j)^{1/4})<\lfloor(n'_j)^{3/4}\rfloor$ for large $j$. Also $n'_j\ge n_*$ for large $j$. Define
\[
W_j=\sum_{k\in I_j}v_{j,k}v_{j,k}^{\mathsf T}.
\]
These real matrices are symmetric positive semidefinite, since
$u^{\mathsf T}W_ju=\sum_{k\in I_j}(u\cdot v_{j,k})^2\ge0$ for every $u\in\R^d$.
Applying \eqref{eq:robust} at scale $n'_j$ and integer $k_j$, and using $I_j\subseteq F_j$, gives
\begin{equation}\label{eq:Wj-energy-bound}
\begin{aligned}
c_{\mathrm{en}}
&\le\sum_{\substack{k\le n'_j\\k\in I_j}}\tn{k_jt_k}^2
=\sum_{\substack{k\le n'_j\\k\in I_j}}(v_{j,k})_{\ell_j}^2\\
&\le\sum_{k\in I_j}|v_{j,k}|^2
=\tr W_j\overset{\eqref{eq:energybound}}{\le}C.
\end{aligned}
\end{equation}
The matrix norm $\|\cdot\|_{\rm op}$ below is induced by the Euclidean norm on $\R^d$. We have
\begin{equation}\label{eq:covariance}
\|W_j\|_{\rm op}=\max\operatorname{spec}(W_j)\le\tr W_j
\overset{\eqref{eq:Wj-energy-bound}}{\le}C.
\end{equation}
Here $\operatorname{spec}(W_j)$ is the set of eigenvalues, all nonnegative by positive semidefiniteness. Finite dimensionality yields a subsequence with $W_j\to\Sigma\ge0$. Continuity of the trace and \eqref{eq:Wj-energy-bound} give
\[
\tr\Sigma=\lim_{j\to\infty}\tr W_j\ge c_{\mathrm{en}}>0,
\]
so $\Sigma\ne0$.

\emph{Separating the independent factor.} Define the centered random vector
$Y_j:\mathcal Z\to\R^d$ by $Y_j=\sum_{k\in I_j}v_{j,k}(Z_k-1/2)$.
Since $b^jt_k-v_{j,k}\in\Z^d$ and $Z_k\in\{0,1\}$, on the torus
\begin{equation}\label{eq:VY}
b^j\sum_{k\le N_j}t_kZ_k=V_j+Y_j\pmod{\Z^d},\qquad
V_j=\sum_{\substack{k\le N_j\\k\notin I_j}}b^jt_kZ_k+\frac12\sum_{k\in I_j}v_{j,k}\pmod{\Z^d}.
\end{equation}
The pair $(V_j,f_j)$ depends only on digits outside $I_j$, hence is independent of $Y_j$. We express their measures as finite Dirac sums to make the convolution algebra explicit. Put $J_j=F_j\setminus I_j$ and $c_j=\frac12\sum_{k\in I_j}v_{j,k}$, and define
\begin{equation}\label{eq:gaussatoms}
\begin{aligned}
\gamma_j&=2^{-|I_j|}\sum_{\zeta\in\{0,1\}^{I_j}}
 \delta_{\rho_d(\sum_{k\in I_j}v_{j,k}(\zeta_k-1/2))},\\
\xi_j&=2^{-|J_j|}\sum_{\zeta\in\{0,1\}^{J_j}}
 \delta_{\rho_d(\sum_{k\in J_j}b^jt_k\zeta_k)},\\
\kappa_j&=(A_*\beta_j)*\xi_j*\delta_{\rho_d(c_j)}.
\end{aligned}
\end{equation}
The first two measures are probabilities; if an index set is empty, its corresponding measure is the point mass at zero. The identity $\gamma_j=(\rho_d\circ Y_j)_*\PP_{\rm dig}$ follows directly from the definition. To identify $\kappa_j$, recall from \eqref{eq:densityprefix} that $f_j(z)=2^{h_j}\nu(E_{z|_{h_j}})$. Since $F_j=\{h_j+1,\ldots,N_j\}$ and $J_j=F_j\setminus I_j$, the map $V_j:\mathcal Z\to\T^d$ in \eqref{eq:VY} can be written as
\[
V_j(z)=\rho_d\left(b^jx_{z|_{h_j}}+\sum_{k\in J_j}b^jt_kz_k+c_j\right),
\qquad x_u=\sum_{k\le h_j}t_ku_k\quad(u\in\{0,1\}^{h_j}).
\]
Under the finite measure $f_j\PP_{\rm dig}$ of mass $m$, a prefix $u$ has mass $2^{-h_j}2^{h_j}\nu(E_u)=\nu(E_u)$. The weight depends only on that prefix, so each binary choice on $J_j$ contributes the factor $2^{-|J_j|}$. Thus
\begin{align*}
(V_j)_*(f_j\PP_{\rm dig})
&=\sum_{u\in\{0,1\}^{h_j}}\nu(E_u)\,2^{-|J_j|}
\sum_{\zeta\in\{0,1\}^{J_j}}
\delta_{\rho_d\left(b^jx_u+\sum_{k\in J_j}b^jt_k\zeta_k+c_j\right)}\\
&=(A_*\beta_j)*\xi_j*\delta_{\rho_d(c_j)}=\kappa_j,
\end{align*}
where the second equality uses \eqref{eq:betaeta}, \eqref{eq:gaussatoms}, and $b^j=Aa^j$. Splitting the independent digits in $F_j$ into $J_j$ and $I_j$, and using $b^jt_k\equiv v_{j,k}\pmod{\Z^d}$ on $I_j$, yields
\[
A_*\eta_j=\xi_j*\delta_{\rho_d(c_j)}*\gamma_j.
\]
Since the homomorphism $A$ preserves convolution, it follows that
\begin{equation}\label{eq:gaussian-convolution-factorization}
\begin{aligned}
A_*(\beta_j*\eta_j)
&=(A_*\beta_j)*(A_*\eta_j)\\
&=(A_*\beta_j)*\xi_j*\delta_{\rho_d(c_j)}*\gamma_j\\
&=\kappa_j*\gamma_j.
\end{aligned}
\end{equation}
\emph{The Gaussian limit on the torus.} For each fixed $\alpha\in\Z^d$, the finite Dirac sum \eqref{eq:gaussatoms} gives
\begin{align*}
\widehat\gamma_j(\alpha)
&=2^{-|I_j|}\sum_{\zeta\in\{0,1\}^{I_j}}
\e\left(\sum_{k\in I_j}(\alpha\cdot v_{j,k})(\zeta_k-\tfrac12)\right)\\
&=2^{-|I_j|}\sum_{\zeta\in\{0,1\}^{I_j}}
\prod_{k\in I_j}\e\left((\alpha\cdot v_{j,k})(\zeta_k-\tfrac12)\right)\\
&=\prod_{k\in I_j}\left(\frac12\sum_{r=0}^1\e\left((\alpha\cdot v_{j,k})(r-\tfrac12)\right)\right)\\
&=\prod_{k\in I_j}\frac{\e(-\tfrac12\alpha\cdot v_{j,k})+\e(\tfrac12\alpha\cdot v_{j,k})}{2}\\
&=\prod_{k\in I_j}\cos(\pi\alpha\cdot v_{j,k}).
\end{align*}
The finite sum factors over the independent binary choices; the last equality uses $(e^{ix}+e^{-ix})/2=\cos x$.
Since $\max_{k\in I_j}|v_{j,k}|\le\delta_j\to0$, all these cosines are positive for sufficiently large $j$. Taylor's formula $\log\cos x=-x^2/2+O(x^4)$ near zero therefore gives
\begin{align*}
\log\widehat\gamma_j(\alpha)
&=-\frac{\pi^2}{2}\sum_{k\in I_j}(\alpha\cdot v_{j,k})^2
 +O\left(\sum_{k\in I_j}|\alpha\cdot v_{j,k}|^4\right)\\
&=-\frac{\pi^2}{2}\alpha^{\mathsf T}W_j\alpha
 +O\left(|\alpha|^4\delta_j^2\sum_{k\in I_j}|v_{j,k}|^2\right)
\longrightarrow-\frac{\pi^2}{2}\alpha^{\mathsf T}\Sigma\alpha.
\end{align*}
The error tends to zero by \eqref{eq:energybound}. Thus
\[
\widehat\gamma_j(\alpha)\longrightarrow
 e^{-\pi^2\alpha^{\mathsf T}\Sigma\alpha/2}
=\widehat\gamma(\alpha),\qquad \gamma=(\rho_d)_*N(0,\Sigma/4),
\]
where the Gaussian coefficient is given by \eqref{eq:gaussianFourier}. All measures have mass one, so Fourier-coefficient convergence implies
\[
\gamma_j\xrightarrow{w^*}\gamma
\]
on $\T^d$, as explained in Section~\ref{sec:prob}.

\emph{Passing to the limit in the factorization.} All $\kappa_j$ have mass $m$. Passing to a further subsequence gives $\kappa_j\xrightarrow{w^*}\kappa$, where $\kappa(\T^d)=m>0$. For each $\alpha\in\Z^d$,
\[
\widehat{\kappa_j*\gamma_j}(\alpha)
=\widehat\kappa_j(\alpha)\widehat\gamma_j(\alpha)
\longrightarrow\widehat\kappa(\alpha)\widehat\gamma(\alpha)
=\widehat{\kappa*\gamma}(\alpha).
\]
Since all these convolutions have mass $m$, the same Fourier-coefficient criterion yields
\[
\kappa_j*\gamma_j\xrightarrow{w^*}\kappa*\gamma.
\]
By the convolution identity \eqref{eq:gaussian-convolution-factorization}, Lemma~\ref{lem:prefix}, and continuity of the fixed homomorphism $A$, the left-hand side also converges to $A_*\sigma$. Hence $A_*\sigma=\kappa*\gamma$, proving \eqref{eq:factor}.
\end{proof}

\begin{proof}[Proof of Theorem~\ref{thm:full}]
Use strong induction over $d\ge1$, assuming the assertion only in dimensions $1\le r<d$. All rank-zero cases are handled directly, so no dimension-zero version is required. Use the arbitrary cluster point $\sigma$ and fixed data chosen in \eqref{eq:fixedlimit}; in particular $\sigma(\T^d)=m>0$. Apply Lemmas~\ref{lem:prefix} and \ref{lem:finiteenergy}. By Lemma~\ref{lem:groups}, $\sigma$ is invariant under $\Gamma^\perp$.

\emph{Case 1: $\rank\Gamma<d$.} If $r=\rank\Gamma>0$, let $A$ have as rows a basis of $\Gamma$. Then $A$ is onto $\T^r$ and $\ker A=\Gamma^\perp$. The sequence $Aa^j$ is quasi-independent: for $\alpha\in\Z^r\setminus\{0\}$, the row $\alpha^{\mathsf T}A$ is a nonzero integer row. For each $\varphi\in C(\T^r)$, testing the original convergence with $\varphi\circ A$ shows that
\[
(m_{Aa^j})_*\nu=A_*((m_{a^j})_*\nu)\xrightarrow{w^*}A_*\sigma.
\]
The induction hypothesis therefore gives $\supp(A_*\sigma)=\T^r$. Invariance under $\ker A$ and \eqref{eq:fibres} then yield
\[
\supp\sigma=A^{-1}(\supp(A_*\sigma))=A^{-1}(\T^r)=\T^d.
\] If $r=0$, invariance under $\Gamma^\perp=\T^d$ and positive mass give the conclusion directly.

\emph{Case 2: $\rank\Gamma=d$.} Lemma~\ref{lem:factor} provides a matrix $\Sigma\ne0$, so Lemma~\ref{lem:gaussgroup} applies to its Gaussian factor. Writing $G=\supp\gamma$, \eqref{eq:convsupp} gives
$\supp(A_*\sigma)=\supp\kappa+G$, a $G$-invariant set. Put $F=\supp\sigma$. Since $\sigma$ is invariant under $\ker A$, \eqref{eq:fibres} gives
\begin{equation}\label{eq:Finvariance}
F=A^{-1}(\supp(A_*\sigma)),\qquad F+A^{-1}G=F.
\end{equation}
Only invariance of the support under $A^{-1}G$ is claimed.
By Lemma~\ref{lem:gaussgroup}, $\Lambda=G^\perp$ has rank $r<d$. If $r>0$, choose its basis matrix $C_\Lambda$ and put $B=C_\Lambda A$. Then $\rank B=r$, $\ker B=A^{-1}G$, and $B:\T^d\to\T^r$ is onto. The sequence $Ba^j$ is quasi-independent by the same integer-row argument. As in the first case, $B_*\sigma$ is a weak-star limit of $(m_{Ba^j})_*\nu$. Induction gives $\supp(B_*\sigma)=\T^r$. Since $F+\ker B=F$ by \eqref{eq:Finvariance}, we may apply \eqref{eq:fibres} with $B$ to conclude
\[
F=B^{-1}(\supp(B_*\sigma))=B^{-1}(\T^r)=\T^d.
\] If $r=0$, Lemma~\ref{lem:gaussgroup} gives $G=\T^d$, and \eqref{eq:Finvariance} yields the conclusion directly. This completes the induction.
\end{proof}

\section{The polar separation}\label{sec:completion}
\begin{proof}[Proof of Theorem~\ref{thm:main}]
For a coefficient sequence in $\mathcal R$, let $K\in H^{(d)}$ and suppose $\mu(K)>0$. Take $\nu=\mu|_K$, and witnesses $(a^j),J$ from Definition~\ref{def:H}. Then
\[
(m_{a^j})_*\nu(J)=\mu(K\cap m_{a^j}^{-1}(J))=0.
\]
The images have common positive mass, so a subsequence converges weak-star to some $\sigma$. The open-set inequality in the Portmanteau theorem gives $\sigma(J)\le\liminf_j(m_{a^j})_*\nu(J)=0$. This contradicts Theorem~\ref{thm:full}, because a positive measure of full support assigns positive mass to every nonempty open set. Thus $\mu\in\HH^\perp$.

The constant $\tau$ in \eqref{eq:tau} is deterministic. Lemmas~\ref{lem:Delta} and \ref{lem:energy} give $\PP_{\rm coef}(\mathcal D_\tau\cap\mathcal R)=1$. Fix any sequence in this intersection. Corollary~\ref{cor:U} gives $E\in U$, while $\mu(E)=1$, so $\mu\notin U^\perp$. The inclusion $\HH\subset U$ gives $U^\perp\subset\HH^\perp$, and our measure makes it strict. Proposition~\ref{prop:synthesis} gives the additional synthesis assertion. The two events concern only the coefficients: no uncountable intersection over measures, dimensions, or sequences is taken.
\end{proof}

\begin{remark}\label{rem:constant}
Every constant choice $\omega_k=c\in[1,2]$ fails the polar conclusion. Indeed, put $K_n=\lfloor n^{1/4}\rfloor$. The pigeonhole principle gives integers $1\le\ell_n\le K_n$ and $a_n$ such that
\[
|a_n-\ell_nQ_n/c|\le K_n^{-1}.
\]
Then $a_n\to\infty$. For $x=\sum_kt_kz_k\in E$, the number $(\ell_nQ_n/c)\sum_{k\le n}t_kz_k$ is an integer. Therefore
\[
\sup_{x\in E}\tn{a_nx}
\le \frac{\sup E}{K_n}+\frac{K_nQ_nL_n}{c}
=O(n^{-1/4})\longrightarrow0.
\]
Thus $E$ is a Dirichlet set: $E\in H^{(1)}\subset U$, whereas $\mu(E)=1$ precludes $\mu\in\HH^\perp$. Moreover $t_1=M_2t_2$, so $b=(1/2,-M_2/2,0,\ldots,0)$ contributes a term tending to one in \eqref{eq:Delta}; the condition $\omega\in\mathcal D_\tau$ fails.
\end{remark}

\section*{Data availability}
No datasets were generated or analysed in this theoretical study.

\section*{Declaration of generative AI and AI-assisted technologies in the manuscript preparation process}
ChatGPT (OpenAI) was used interactively to explore proof strategies, develop intermediate arguments and estimates, check calculations, and assist with literature searches. These interactions informed the development of the proofs presented in this article; AI output is not used as a substitute for a mathematical proof or as an independent verification of a claim.

During the preparation of this manuscript, the author used ChatGPT (OpenAI) to assist with English translation, wording, organization, and LaTeX preparation. The author independently verified all mathematical arguments, reviewed and revised the AI-assisted material, and takes full responsibility for the content of the article.

% Funding and competing-interest declarations require the author's factual input;
% no absence of funding or competing interests is asserted in this draft.
% BEGIN inlined references.tex

% END inlined references.tex
\end{document}